\documentclass[12pt]{article}
\usepackage{amsmath,amsfonts,amssymb,amsthm}
\input amssym.def
\oddsidemargin
0.5cm

\begin{document}
\def \Z{\Bbb Z}
\def \C{\Bbb C}
\def \R{\Bbb R}
\def \Q{\Bbb Q}
\def \N{\Bbb N}
\def \F{\Bbb F}

\def \A{{\mathcal{A}}}
\def \D{{\mathcal{D}}}
\def \E{{\mathcal{E}}}
\def \L{\mathcal{L}}
\def \S{{\mathcal{S}}}
\def \wt{{\rm wt}}
\def \tr{{\rm tr}}
\def \span{{\rm span}}
\def \Res{{\rm Res}}
\def \Der{{\rm Der}}
\def \End{{\rm End}}
\def \Ind {{\rm Ind}}
\def \Irr {{\rm Irr}}
\def \Aut{{\rm Aut}}
\def \GL{{\rm GL}}
\def \Hom{{\rm Hom}}
\def \mod{{\rm mod}}
\def \ann{{\rm Ann}}
\def \ad{{\rm ad}}
\def \rank{{\rm rank}\;}
\def \<{\langle}
\def \>{\rangle}

\def \g{{\frak{g}}}
\def \h{{\hbar}}
\def \k{{\frak{k}}}
\def \sl{{\frak{sl}}}
\def \gl{{\frak{gl}}}

\def \be{\begin{equation}\label}
\def \ee{\end{equation}}
\def \bex{\begin{example}\label}
\def \eex{\end{example}}
\def \bl{\begin{lem}\label}
\def \el{\end{lem}}
\def \bt{\begin{thm}\label}
\def \et{\end{thm}}
\def \bp{\begin{prop}\label}
\def \ep{\end{prop}}
\def \br{\begin{rem}\label}
\def \er{\end{rem}}
\def \bc{\begin{coro}\label}
\def \ec{\end{coro}}
\def \bd{\begin{de}\label}
\def \ed{\end{de}}

\newcommand{\m}{\bf m}
\newcommand{\n}{\bf n}
\newcommand{\nno}{\nonumber}
\newcommand{\nord}{\mbox{\scriptsize ${\circ\atop\circ}$}}
\newtheorem{thm}{Theorem}[section]
\newtheorem{prop}[thm]{Proposition}
\newtheorem{coro}[thm]{Corollary}
\newtheorem{conj}[thm]{Conjecture}
\newtheorem{example}[thm]{Example}
\newtheorem{lem}[thm]{Lemma}
\newtheorem{rem}[thm]{Remark}
\newtheorem{de}[thm]{Definition}
\newtheorem{hy}[thm]{Hypothesis}
\makeatletter
\@addtoreset{equation}{section}
\def\theequation{\thesection.\arabic{equation}}
\makeatother
\makeatletter

\begin{center}
{\Large \bf Quantum vertex $\F((t))$-algebras and their modules }
\end{center}

\begin{center}
{Haisheng Li\footnote{Partially supported
by NSF grant DMS-0600189}\\
Department of Mathematical Sciences\\
Rutgers University, Camden, NJ
08102}
\end{center}

\begin{abstract}
This is a paper in a series to study vertex algebra-like structures
arising from various algebras including quantum affine algebras and
Yangians. In this paper, we develop a theory of what we call (weak)
quantum vertex $\F((t))$-algebras with $\F$ a field of
characteristic zero and $t$ a formal variable, and we give a
conceptual construction of (weak) quantum vertex $\F((t))$-algebras
and their modules. As an application, we associate weak quantum
vertex $\F((t))$-algebras to quantum affine algebras, providing a
solution to a problem posed by Frenkel and Jing.  We also explicitly
construct an example of quantum vertex $\F((t))$-algebras from a
certain quantum $\beta\gamma$-system.
\end{abstract}

\section{Introduction}
In the earliest days of vertex (operator) algebra theory, Lie
algebras had played an important role. In particular, an important
family of vertex operator algebras (see \cite{flm}, \cite{fz},
\cite{dl}) was associated to untwisted affine Lie algebras. A
fundamental problem, posed in \cite{fj} (see also \cite{efk}), has
been to establish a suitable theory of quantum vertex algebras
 so that quantum vertex algebras can be
canonically associated to quantum affine algebras in the same (or
similar) way that vertex operator algebras are associated to affine
Lie algebras. In the past, several theories of quantum vertex
algebras have been studied (\cite{efr}, \cite{ek}, \cite{b-qva},
\cite{li-qva1}, \cite{li-qva2}, \cite{ab}), however this particular
problem is still to be solved.

This is a paper in a series, starting with \cite{li-qva1}, to study
vertex algebra-like structures arising from various algebras such as
quantum affine algebras and Yangians, with an ultimate goal to solve
the aforementioned problem. In the present paper, we develop a
theory of (weak) quantum vertex $\F((t))$-algebras with $\F$ a field
of characteristic zero and $t$ a formal variable, and we establish a
general construction of weak quantum vertex $\F((t))$-algebras and
their modules. As an application we associate weak quantum vertex
$\F((t))$-algebras canonically to quantum affine algebras, providing
a desired solution to the very problem.

The notion of weak quantum vertex $\F((t))$-algebra in a certain way
generalizes the notion of weak quantum vertex algebra, which was
introduced and studied previously in this series (see
\cite{li-qva1}, \cite{li-qva2}). A rough description of all these
``quantum vertex algebras'' is that they are various generalizations
of ordinary vertex algebras where the locality, namely weak
commutativity, is replaced by a braided locality, while the weak
associativity is retained. A {\em weak quantum vertex
$\F((t))$-algebra} is defined to be an $\F((t))$-module $V$,
equipped with an $\F$-linear map
$$Y(\cdot,x): V\rightarrow \Hom_{\F((t))} (V,V((x)))$$
and equipped with a distinguished vector ${\bf 1}\in V$, satisfying
the conditions that
$$Y(f(t)v,x)=f(t+x)Y(v,x)\ \ \ \mbox{
for }f(t)\in \F((t)),\; v\in V,$$
$$Y({\bf 1},x)v=v,\ \ Y(v,x){\bf 1}\in V[[x]]\ \ \mbox{and }\ \lim_{x\rightarrow
0}Y(v,x){\bf 1}=v\ \ \mbox{ for }v\in V,$$ and that for $u,v\in V$,
there exist (finitely many)
$$u^{(i)},\; v^{(i)}\in V,\; f_{i}(x_{1},x_{2})\in
\F_{*}(x_{1},x_{2})\ \ (i=1,\dots,r)$$ such that
\begin{eqnarray*}
& &x_{0}^{-1}\delta\left(\frac{x_{1}-x_{2}}{x_{0}}\right)
Y(u,x_{1})Y(v,x_{2})\\
& &\ \ \ \ -x_{0}^{-1}\delta\left(\frac{x_{2}-x_{1}}{-x_{0}}\right)
\sum_{i=1}^{r}\iota_{t,x_{2},x_{1}}(f_{i}(x_{1}+t,x_{2}+t))
Y(v^{(i)},x_{2})Y(u^{(i)},x_{1})\\
 &=&x_{2}^{-1}\delta\left(\frac{x_{1}-x_{0}}{x_{2}}\right)
Y(Y(u,x_{0})v,x_{2}).
\end{eqnarray*}
(See Section 2 for the definitions of $\F_{*}(x_{1},x_{2})$ and
$\iota_{t, x_{2},x_{1}}$.) Furthermore, a quantum vertex
$\F((t))$-algebra is a weak quantum vertex $\F((t))$-algebra
equipped with a unitary quantum Yang-Baxter operator on $V$ with two
(independent) spectral parameters, which describes the braiding and
satisfies some other conditions.

In \cite{ek}, Etingof and Kazhdan developed a fundamental theory of
quantum vertex operator algebras in the sense of formal deformation.
The notion of (weak) quantum vertex $\F((t))$-algebra as well as
that of (weak) quantum vertex algebra (see \cite{li-qva1},
\cite{li-qva2}) largely reflects Etingof-Kazhdan's notion of quantum
vertex operator algebra, however there are essential differences. As
the map $Y(\cdot,x)$ for a weak quantum vertex $\F((t))$-algebra is
{\em not} $\F((t))$-linear (where linearity is deformed), the formal
variable $t$ is not a deformation parameter, unlike the formal
variable $\hbar$ in Etingof-Kazhdan's theory. On the other hand, the
braiding operator in Etingof-Kazhdan's theory is a rational quantum
Yang-Baxter operator (with one parameter), whereas the braiding
operator here is more general with two parameters.

The theory of quantum vertex $\F((t))$-algebras is also
significantly different from Anguelova and Bergvelt's theory of
$H_{D}$-quantum vertex algebras (see \cite{ab}). The notion of
$H_{D}$-quantum vertex algebra generalizes Etingof-Kazhdan's notion
of braided vertex operator algebra (see \cite{ek}) in certain
directions. In particular, the underlying space of an
$H_{D}$-quantum vertex algebra is a topologically free
$\F[[t]]$-module and the vertex operator map $Y(\cdot,x)$ is
$\F[[t]]$-linear, where the variable $t$ plays the same role as
$\hbar$ does in \cite{ek}. We note that weak quantum vertex
$\F((t))$-algebras satisfy the same associativity for ordinary
vertex algebras. Unlike (weak) quantum vertex $\F((t))$-algebras,
general $H_{D}$-quantum vertex algebras {\em do not} satisfy the
associativity for ordinary vertex algebras (though they do satisfy a
braided associativity).

The theory of (weak) quantum vertex $\F((t))$-algebras is deeply
rooted in \cite{li-qva1}. To better state the results of the present
paper we review a conceptual result obtained therein. Let $W$ be an
{\em arbitrary} vector space and let $\E(W)$ denote the space $\Hom
(W,W((x)))$ alternatively. The essential idea is to study the
algebraic structures generated by various types of subsets of
$\E(W)$. The most general type consists of what we called quasi
compatible subsets, where a subset $U$ of $\E(W)$ is {\em quasi
compatible} if for any finite sequence $a^{(1)}(x),\dots,a^{(r)}(x)$
in $U$, there exists a nonzero polynomial $p(x,y)$ such that
\begin{eqnarray*}
\left(\prod_{1\le i<j\le r}p(x_{i},x_{j})\right) a_{1}(x_{1})\cdots
a_{r}(x_{r})\in \Hom (W,W((x_{1},\dots,x_{r}))).
\end{eqnarray*}
Furthermore, the notion of {\em compatible subset} is defined by
assuming that the nonzero polynomial $p(x,y)$ is of the form
$(x-y)^{k}$ with $k\in \N$. It was proved therein that any (quasi)
compatible subset $U$ of $\E(W)$ generates what we called a nonlocal
vertex algebra $\<U\>$ in a certain canonical way with $W$ as a
(quasi) module in a certain sense. (A nonlocal vertex algebra is the
same as a weak $G_{1}$-vertex algebra in the sense of \cite{li-g1}
and is also essentially the same as a field algebra in the sense of
\cite{bk}.) In contrast with that vertex algebras are analogs of
commutative and associative algebras, nonlocal vertex algebras are
analogs of noncommutative associative algebras. It follows from this
general result that nonlocal vertex algebras can be associated to a
wide variety of algebras including quantum affine algebras.

In the present paper, based on \cite{li-qva1}, as one of our main
results we prove that for any quasi compatible subset $U$ of
$\E(W)$, the $\F((x))$-span $\F((x))\<U\>$ is what we call a
nonlocal vertex $\F((t))$-algebra. (Note that $\E(W)$ is naturally
an $\F((x))$-module.) The notion of nonlocal vertex
$\F((t))$-algebra is a counterpart of the notion of nonlocal vertex
algebra, where a nonlocal vertex $\F((t))$-algebra $V$ is a nonlocal
vertex algebra over $\F$ and an $\F((t))$-module such that
$$Y(f(t)u,x)g(t)v=f(t+x)g(t)Y(u,x)v\ \ \ \mbox{
for }f(t),g(t)\in \F((t)),\; u,v\in V.$$ Furthermore, to deal with
quantum affine algebras, we study what we call quasi
$\S(x_{1},x_{2})$-local subsets of $\E(W)$. A subset $U$ of $\E(W)$
is said to be {\em quasi $\S(x_{1},x_{2})$-local}  if for any
$a(x),b(x)\in U$, there exist (finitely many)
$$u^{(i)}(x),v^{(i)}(x)\in U,\; f_{i}(x_{1},x_{2})\in \F_{*}(x_{1},x_{2})
\ \ (i=1,\dots,r)$$ such that
$$p(x_{1},x_{2})a(x_{1})b(x_{2})
=p(x_{1},x_{2})\sum_{i=1}^{r}\iota_{x_{2},x_{1}}(f_{i}(x_{1},x_{2}))
u^{(i)}(x_{2})v^{(i)}(x_{1})$$ for some nonzero polynomial
$p(x_{1},x_{2})$. We note that quasi $\S(x_{1},x_{2})$-local subsets
are quasi compatible. Our key result is that for any quasi
$\S(x_{1},x_{2})$-local subset $U$ of $\E(W)$, $\F((x))\<U\>$ is a
weak quantum vertex $\F((t))$-algebra.

The theory of (weak) quantum vertex $\F((t))$-algebras runs largely
parallel to that of (weak) quantum vertex algebras. To construct
quantum vertex $\F((t))$-algebras from weak quantum vertex
$\F((t))$-algebras, we extend Etingof-Kazhdan's notion of
non-degeneracy for nonlocal vertex $\F((t))$-algebras and we prove
that every non-degenerate weak quantum $\F((t))$-algebra has a
(unique) canonical quantum vertex $\F((t))$-algebra structure, just
as with weak quantum vertex algebras in \cite{li-qva2} (see also
\cite{ek}).  We furthermore establish certain general non-degeneracy
results, analogous to those obtained in \cite{li-qva2}.

We note that this theory of quantum vertex $\F((t))$-algebras has a
great generality and our conceptual result is applicable to many
better known quantum algebras, particularly including quantum affine
algebras. Take $W$ to be a highest weight module for a quantum
affine algebra and set $\F=\C$. We show that the generating
functions of the generators in Drinfeld's realization form a quasi
$\S(x_{1},x_{2})$-local subset $U$ of $\E(W)$. Therefore we have a
weak quantum vertex $\C((t))$-algebra $\C((x))\<U\>$ with $W$ as a
canonical quasi module. To a certain extent, this solves the
aforementioned problem, though we yet have to show that this weak
quantum vertex $\C((t))$-algebra is a quantum vertex
$\C((t))$-algebra, or sufficiently to show that it is
non-degenerate.

In the theory of (weak) quantum vertex $\F((t))$-algebras, an
important issue is about notions of module. Notice that for a quasi
$\S(x_{1},x_{2})$-local subset $U$ of $\E(W)$ with $W$ a vector
space as before, the weak quantum vertex $\F((t))$-algebra
$\F((x))\<U\>$ has the natural module $W$ (a vector space over $\F$)
and the adjoint module $\F((x))\<U\>$ (a vector space over
$\F((t))$), which are significantly different. This leads to us to
two categories of modules for weak quantum vertex
$\F((t))$-algebras.

This paper is organized as follows: In Section 2, we study notions
of nonlocal vertex $\F((t))$-algebra and weak quantum vertex
$\F((t))$-algebra. In Section 3, we study notions of quantum vertex
$\F((t))$-algebra and non-degeneracy. In Section 4, we give a
conceptual construction of nonlocal vertex $\F((t))$-algebras and
weak quantum vertex $\F((t))$-algebras. In Section 5, we present two
existence theorems. In Section 6, we associate weak quantum vertex
$\C((t))$-algebras to quantum affine algebras and we construct a
quantum vertex $\C((t))$-algebra {}from a certain quantum
$\beta\gamma$-system.

\section{Nonlocal vertex $\F((t))$-algebras and weak quantum vertex $\F((t))$-algebras}
In this section, we define notions of nonlocal vertex
$\F((t))$-algebra and weak quantum vertex $\F((t))$-algebra, and we
study what we call type zero modules and type one modules for
nonlocal vertex $\F((t))$-algebras. We also present some basic
axiomatic results.

We begin by fixing some basic notations. In addition to the standard
usage of symbols $\Z$ and $\C$,  we use $\N$ for the set of
nonnegative integers. We shall use the standard formal variable
notations and conventions as in \cite{flm} and \cite{fhl} (cf.
\cite{ll}). Letters such as $t, x,y, z, x_{0},x_{1}, x_{2},\dots$
stand for mutually commuting independent formal variables.  We shall
be working on a scalar field $\F$ of characteristic zero, where
typical examples of $\F$ are $\C$ and the field $\C((t))$ of lower
truncated formal Laurent series in $t$. Denote by
$\F((x_{1},\dots,x_{r}))$ the algebra of formal Laurent series which
are globally truncated with respect to all the variables. By
$\F_{*}(x_{1},x_{2},\dots,x_{r})$ we denote the extension of the
algebra $\F[[x_{1},x_{2},\dots,x_{r}]]$ of formal nonnegative power
series by joining the inverses of nonzero polynomials.

We recall the iota maps from \cite{li-qva1}, which will be used
extensively. For any permutation $(i_{1},i_{2},\dots,i_{r})$ on $\{
1,\dots,r\}$,
\begin{eqnarray}
\iota_{x_{i_{1}},\dots,x_{i_{r}}}: \F_{*}(x_{1},x_{2},\dots,x_{r})
\rightarrow \F((x_{i_{1}}))\cdots ((x_{i_{r}}))
\end{eqnarray}
denotes the unique algebra embedding that extends the identity
endomorphism of $\F[[x_{1},\dots,x_{r}]]$ (cf. \cite{fhl}). Note
that both $\F_{*}(x_{1},\dots,x_{r})$ and $\F((x_{i_{1}}))\cdots
((x_{i_{r}}))$ contain $\F((x_{1},\dots,x_{r}))$ as a subalgebra.
The map $\iota_{x_{i_{1}},\dots,x_{i_{r}}}$ preserves
$\F((x_{1},\dots,x_{r}))$ element-wise and is
$\F((x_{1},\dots,x_{r}))$-linear.

We recall the notion of nonlocal vertex algebra (\cite{li-g1},
\cite{li-qva1}; see also \cite{kac}, \cite{bk}), which is essential
to this paper.

\bd{dnonlocalva} {\em  A {\em nonlocal vertex algebra} over $\F$ is
a vector space $V$, equipped with a linear map
\begin{eqnarray*}
Y(\cdot,x): &&V\rightarrow \Hom (V,V((x)))\subset (\End V)[[x,x^{-1}]],\\
&&v\mapsto Y(v,x)=\sum_{n\in \Z}v_{n}x^{-n-1}\ \ (\mbox{with
}v_{n}\in \End V)
\end{eqnarray*}
and a distinguished vector ${\bf 1}\in V$, satisfying the conditions
that
$$ Y({\bf 1},x)v=v,\ \ Y(v,x){\bf 1}\in V[[x]]
\ \ \mbox{and }\ \lim_{x\rightarrow 0}Y(v,x){\bf 1}=v
\ \ \mbox{ for }v\in V,$$ and that for $u,v,w\in V$, there exists a
nonnegative integer $l$ such that
$$(x_{0}+x_{2})^{l}Y(u,x_{0}+x_{2})Y(v,x_{2})w
=(x_{0}+x_{2})^{l}Y(Y(u,x_{0})v,x_{2})w$$ (the {\em weak
associativity}).} \ed

The following two notions can be found either in \cite{li-g1} or
\cite{li-qva1}:

\bd{dmodule} {\em Let $V$ be a nonlocal vertex algebra. A {\em
$V$-module} is a vector space $W$ equipped with a linear map
\begin{eqnarray*}
Y_{W}(\cdot,x):&& V\rightarrow \Hom (W,W((x)))\subset (\End
W)[[x,x^{-1}]],\\
&&v\mapsto Y_{W}(v,x),
\end{eqnarray*}
satisfying the conditions that
$$Y_{W}({\bf 1},x)=1_{W}\ \ \mbox{
(the identity operator on $W$)}$$
 and that for $u,v\in V,\; w\in W$,
there exists a nonnegative integer $l$ such that
$$(x_{0}+x_{2})^{l}Y_{W}(u,x_{0}+x_{2})Y_{W}(v,x_{2})w
=(x_{0}+x_{2})^{l}Y_{W}(Y(u,x_{0})v,x_{2})w.$$ The notion of {\em
quasi $V$-module} is defined as above with the last condition
replaced by a weaker condition that for $u,v\in V,\; w\in W$, there
exists a nonzero polynomial $p(x_{1},x_{2})\in \F[x_{1},x_{2}]$ such
that
\begin{eqnarray}
p(x_{0}+x_{2},x_{2})Y_{W}(u,x_{0}+x_{2})Y_{W}(v,x_{2})w
=p(x_{0}+x_{2},x_{2})Y_{W}(Y(u,x_{0})v,x_{2})w.
\end{eqnarray}}
\ed

The following follows immediately from \cite{ltw} (Lemma 2.9):

\bp{prepresentation-F} Let $V$ be a nonlocal vertex algebra. In the
definition of a $V$-module, in the presence of other axioms weak
associativity can be equivalently replaced by the condition that for
$u,v\in V$, there exists $k\in \N$ such that
\begin{eqnarray*}
(x_{1}-x_{2})^{k}Y_{W}(u,x_{1})Y_{W}(v,x_{2})\in
\Hom(W,W((x_{1},x_{2}))),
\end{eqnarray*}
$$x_{0}^{k}Y_{W}(Y(u,x_{0})v,x)
=\left((x_{1}-x_{2})^{k}Y_{W}(u,x_{1})Y_{W}(v,x_{2})\right)|_{x_{1}=x_{2}+x_{0}}.$$
\ep

\br{rsubstitution} {\em  For $A(x_{1},x_{2})\in \Hom
(W,W((x_{1}))((x_{2})))$ with $W$ a vector space over $\F$, we have
been using the convention
$$A(x_{1},x_{2})|_{x_{1}=x_{0}+x_{2}}=A(x_{0}+x_{2},x_{2})
=\iota_{x_{0},x_{2}}A(x_{0}+x_{2},x_{2}).$$ Note that the
substitutions $A(x_{2}+x_{0},x_{2})$, $A(x_{1},x_{1}+x_{0})$ and
$A(x_{1},x_{0}+x_{1})$ do {\em not} exist in general. On the other
hand, for $E(x_{1},x_{2})\in \Hom (W,W((x_{1},x_{2})))$, all the
substitutions $E(x_{0}+x_{2},x_{2})$, $E(x_{2}+x_{0},x_{2})$,
$E(x_{1},x_{1}-x_{0})$, and $E(x_{1},-x_{0}+x_{1})$ exist, and we
have
\begin{eqnarray}
\left(E(x_{1},x_{2})|_{x_{2}=x_{1}-x_{0}}\right)|_{x_{1}=x_{2}+x_{0}}
=E(x_{1},x_{2})|_{x_{1}=x_{2}+x_{0}}.
\end{eqnarray}}
\er

Let $\F$ be a field of characteristic zero as before and let $t$ be
a formal variable. Notice that as $\F((t))$ is a field containing
$\F$ as a subfield, every $\F((t))$-module is naturally a vector
space over $\F$.

\bd{dt-algebra} {\em A {\em nonlocal vertex $\F((t))$-algebra} is a
nonlocal vertex algebra $V$ over $\F$, equipped with an
$\F((t))$-module structure, such that
\begin{eqnarray}
Y(f(t)u,x)(g(t)v)=f(t+x)g(t)Y(u,x)v
\end{eqnarray}
for $f(t),g(t)\in \F((t)),\; u, v\in V$, where it is understood that
$$f(t+x)=e^{x\frac{d}{dt}}f(t)\in \F((t))[[x]].$$ A {\em homomorphism} of
nonlocal vertex $\F((t))$-algebras is a homomorphism of nonlocal
vertex algebras over $\F$, which is also $\F((t))$-linear.} \ed

\bd{dmodule} {\em Let $V$ be a nonlocal vertex $\F((t))$-algebra. A
{\em  $V$-module of type zero} is a module $(W,Y_{W})$ for $V$
viewed as a nonlocal vertex algebra over $\F$, satisfying the
condition that
\begin{eqnarray}
Y_{W}(f(t)v,x)w=f(x)Y_{W}(v,x)w \ \ \ \mbox{ for }f(t)\in \F((t)),\;
v\in V,\; w\in W.
\end{eqnarray}
We define a notion of {\em quasi $V$-module of type zero} in the
obvious way---with the word ``module'' replaced by ``quasi module''
in the two places.} \ed

The following immediately follows from the corresponding results for
nonlocal vertex algebras (see \cite{li-qva1}):

\bl{ld-property} Let $V$ be a nonlocal vertex $\F((t))$-algebra and
let $\D$ be the $\F$-linear operator on $V$ defined by $\D
v=v_{-2}{\bf 1}$ for $v\in V$. Then
\begin{eqnarray}
& &Y(v,x){\bf 1}=e^{x\D}v,\label{ecreation}\\
& &[\D,Y(v,x)]=Y(\D v,x)={d\over dx}Y(v,x),\label{ed-bracket}\\
& & e^{x\D} (f(t)v)=f(t+x)e^{x\D} v\ \ \ \ \ \ \mbox{ for }f(t)\in
\F((t)),\; v\in V.
\end{eqnarray}
Furthermore, for any type-zero quasi $V$-module $(W,Y_{W})$ we have
\begin{eqnarray}
Y_{W}(\D v,x)=\frac{d}{dx}Y_{W}(v,x)\ \ \ \mbox{ for }v\in V.
\end{eqnarray}
\el

Note that as ${\partial\over
\partial x_{1}}+{\partial\over \partial x_{2}}$ is a derivation of
$\F[[x_{1},x_{2}]]$, $e^{t(\partial/\partial x_{1}+\partial/\partial
x_{2})}$ is an algebra embedding of $\F[[x_{1},x_{2}]]$ into
$\F[[t,x_{1},x_{2}]]$ with
$$e^{x\left({\partial\over \partial x_{1}}+{\partial\over \partial x_{2}}\right)}
\F[x_{1},x_{2}] \subset \F[t,x_{1},x_{2}].$$ Consequently, this
gives rise to an algebra embedding of $\F_{*}(x_{1},x_{2})$ into
$\F_{*}(t,x_{1},x_{2})$, where for $f(x_{1},x_{2})\in
\F_{*}(x_{1},x_{2})$,
$$e^{t\left({\partial\over \partial x_{1}}+{\partial\over \partial x_{2}}\right)}
f(x_{1},x_{2})=f(x_{1}+t,x_{2}+t)\in \F_{*}(t,x_{1},x_{2}).$$

We now define the main object of this paper.

\bd{dtqva} {\em A {\em weak quantum vertex $\F((t))$-algebra} is a
nonlocal vertex $\F((t))$-algebra $V$, satisfying the condition that
for any $u,v\in V$, there exist
$$u^{(i)},v^{(i)}\in V,\; f_{i}(x_{1},x_{2})\in \F_{*}(x_{1},x_{2})
\ \ \ (i=1,\dots,r)$$
 such that
\begin{eqnarray}\label{et-jacobi-wqva}
& &x_{0}^{-1}\delta\left(\frac{x_{1}-x_{2}}{x_{0}}\right)
Y(u,x_{1})Y(v,x_{2})\nonumber\\
& &\ \ \ \ - x_{0}^{-1}\delta\left(\frac{x_{2}-x_{1}}{-x_{0}}\right)
\sum_{i=1}^{r}\iota_{t,x_{2},x_{1}}(f_{i}(t+x_{1},t+x_{2}))
Y(v^{(i)},x_{2})Y(u^{(i)},x_{1})\nonumber\\
&=&x_{2}^{-1}\delta\left(\frac{x_{1}-x_{0}}{x_{2}}\right)
Y(Y(u,x_{0})v,x_{2})
\end{eqnarray}
(the {\em $\S_{t}(x_{1},x_{2})$-Jacobi identity}). } \ed

In the following we study certain axiomatic aspects. For convenience
we recall from \cite{li-g1} the following result (cf. \cite{fhl}):

\bl{lformal-jacobi-old} Let $W$ be a vector space over $\F$ and let
\begin{eqnarray*}
A(x_{1},x_{2})\in W((x_{1}))((x_{2})),\  B(x_{1},x_{2})\in
W((x_{2}))((x_{1})),\
%&&\hspace{3cm}
C(x_{0},x_{2})\in W((x_{2}))((x_{0})).
\end{eqnarray*}
Then
\begin{eqnarray*}\label{ejacobi-abc}
&&x_{0}^{-1}\delta\left(\frac{x_{1}-x_{2}}{x_{0}}\right)
A(x_{1},x_{2})-
x_{0}^{-1}\delta\left(\frac{x_{2}-x_{1}}{-x_{0}}\right)
B(x_{1},x_{2})\\
&&\hspace{2cm}=x_{2}^{-1}\delta\left(\frac{x_{1}-x_{0}}{x_{2}}\right)
C(x_{0},x_{2})
\end{eqnarray*}
if and only if there exist nonnegative integers $k$ and $l$ such
that
\begin{eqnarray*}
&&(x_{1}-x_{2})^{k}A(x_{1},x_{2})=(x_{1}-x_{2})^{k}B(x_{1},x_{2}),\\
&&(x_{0}+x_{2})^{l}C(x_{0},x_{2})
=(x_{0}+x_{2})^{l}A(x_{0}+x_{2},x_{2}).
\end{eqnarray*}
\el

As weak associativity holds for every nonlocal vertex algebra, in
view of Lemma \ref{lformal-jacobi-old} we immediately have:

\bp{ptlocality} In Definition \ref{dtqva}, the
$\S_{t}(x_{1},x_{2})$-Jacobi identity axiom in the presence of other
axioms can be equivalently replaced by {\em
$\S_{t}(x_{1},x_{2})$-locality:} For $u,v\in V$, there exist
$$u^{(i)},v^{(i)}\in V,\; f_{i}(x_{1},x_{2})\in \F_{*}(x_{1},x_{2})
\ \ \ (i=1,\dots,r)$$  such that
\begin{eqnarray}\label{et-local}
& &(x_{1}-x_{2})^{k}Y(u,x_{1})Y(v,x_{2})\nonumber\\
&=&\sum_{i=1}^{r}(x_{1}-x_{2})^{k}
\iota_{t,x_{2},x_{1}}(f_{i}(x_{1}+t,x_{2}+t))
Y(v^{(i)},x_{2})Y(u^{(i)},x_{1})
\end{eqnarray}
for some nonnegative integer $k$ depending on $u$ and $v$. \ep

We also have:

\bp{pskewsymmetry} Let $V$ be a nonlocal vertex $\F((t))$-algebra
and let
$$u,\ v,\ u^{(i)},\ v^{(i)}\in V,\; f_{i}(x_{1},x_{2})\in \F_{*}(x_{1},x_{2})
\ \ (i=1,\dots,r).$$ Then (\ref{et-local}) holds for some
nonnegative integer $k$ if and only if
\begin{eqnarray}\label{eskew-symmetry-1}
Y(u,x)v&=&\sum_{i=1}^{r}\iota_{t,x}(f_{i}(x+t,t))
e^{x\D}Y(v^{(i)},-x)u^{(i)}.
\end{eqnarray}
\ep

\begin{proof} We follow the proof of the analogous assertion
for ordinary vertex algebras in \cite{ll}. Assume that
(\ref{et-local}) holds for some $k\in \N$. We can choose $k$ so
large that
$$x^{k}Y(v^{(i)},x)u^{(i)}\in V[[x]]\ \ \ \mbox{ for }i=1,\dots,r.$$
Let $p(x_{1},x_{2})\in \F[x_{1},x_{2}]$ be a nonzero polynomial such
that
$$p(x_{1},x_{2})f_{i}(x_{1},x_{2})\in \F[[x_{1},x_{2}]]
\ \ \ \mbox{ for }i=1,\dots,r.$$
 Using (\ref{et-local}) and
the $\D$-properties in Lemma \ref{ld-property} we have
\begin{eqnarray*}
& &p(x_{1}+t,x_{2}+t)(x_{1}-x_{2})^{k}Y(u,x_{1})Y(v,x_{2}){\bf 1}\nonumber\\
&=&\sum_{i=1}^{r}(x_{1}-x_{2})^{k}p(x_{1}+t,x_{2}+t)
\iota_{t,x_{2},x_{1}}(f_{i}(x_{1}+t,x_{2}+t))
Y(v^{(i)},x_{2})Y(u^{(i)},x_{1}){\bf 1}\\
&=&\sum_{i=1}^{r}(x_{1}-x_{2})^{k} (pf_{i})(x_{1}+t,x_{2}+t)
Y(v^{(i)},x_{2})e^{x_{1}\D}u^{(i)}\\
&=&\sum_{i=1}^{r}(x_{1}-x_{2})^{k} (pf_{i})(x_{1}+t,x_{2}+t)
e^{x_{1}\D}Y(v^{(i)},x_{2}-x_{1})u^{(i)}\\
&=&\sum_{i=1}^{r}(pf_{i})(x_{1}+t,x_{2}+t)
e^{x_{1}\D}[(x_{1}-x_{2})^{k}Y(v^{(i)},-x_{1}+x_{2})u^{(i)}].
\end{eqnarray*}
Notice that it is safe now to set $x_{2}=0$. By doing so we get
\begin{eqnarray*}
& &p(x_{1}+t,t)x_{1}^{k}Y(u,x_{1})v\\
&=&\sum_{i=1}^{r}x_{1}^{k}(pf_{i})(x_{1}+t,t)
e^{x_{1}\D}Y(v^{(i)},-x_{1})u^{(i)}\\
&=&\sum_{i=1}^{r}x_{1}^{k}p(x_{1}+t,t)\iota_{t,x_{1}}(f_{i}(x_{1}+t,t))
e^{x_{1}\D}Y(v^{(i)},-x_{1})u^{(i)}.
\end{eqnarray*}
By cancellation (namely, multiplying both sides by
$\iota_{t,x_{1}}(x_{1}^{-k}/p(x_{1}+t,t))$) we obtain
\begin{eqnarray*}
Y(u,x_{1})v &=&\sum_{i=1}^{r}\iota_{t,x_{1}}(f_{i}(x_{1}+t,t))
e^{x_{1}\D}Y(v^{(i)},-x_{1})u^{(i)}.
\end{eqnarray*}

On the other hand, assume that this skew-symmetry relation holds. By
Proposition \ref{prepresentation-F}, there exists $k\in \N$ such
that
\begin{eqnarray*}
& &(x_{1}-x_{2})^{k}Y(u,x_{1})Y(v,x_{2})\in \Hom (V,V((x_{1},x_{2}))),\\
& &x_{0}^{k}Y(Y(u,x_{0})v,x_{2})
=\left((x_{1}-x_{2})^{k}Y(u,x_{1})Y(v,x_{2})\right)|_{x_{1}=x_{2}+x_{0}};\\
& &(x_{1}-x_{2})^{k}Y(v^{(i)},x_{2})Y(u^{(i)},x_{1})\in \Hom (V,V((x_{1},x_{2}))),\\
& &x_{0}^{k}Y(Y(v^{(i)},-x_{0})u^{(i)},x_{1})
=\left((x_{1}-x_{2})^{k}Y(v^{(i)},x_{2})Y(u^{(i)},x_{1})w\right)|_{x_{2}=x_{1}-x_{0}},
\end{eqnarray*}
and such that
$$\iota_{t,x_{1},x_{2}}(x_{1}-x_{2})^{k}f_{i}(t+x_{1},t+x_{2})
=\iota_{t,x_{2},x_{1}}(x_{1}-x_{2})^{k}f_{i}(t+x_{1},t+x_{2}),$$
lying in $\F((t))[[x_{1},x_{2}]]$ for $i=1,\dots,r$ (recall Lemma
\ref{lsubcomm}). Set
$$E(x_{1},x_{2})=\sum_{i=1}^{r}\iota_{t,x_{2},x_{1}}(f_{i}(t+x_{1},t+x_{2}))
(x_{1}-x_{2})^{2k}Y(v^{(i)},x_{2})Y(u^{(i)},x_{1}).$$ Then
$$E(x_{1},x_{2})\in \Hom (V,V((x_{1},x_{2}))).$$
 Using the skew-symmetry
relation and the basic $\D$-properties we get
\begin{eqnarray*}
& &Y(Y(u,x_{0})v,x_{2})\\
&=&\sum_{i=1}^{r}Y\left(\iota_{t,x_{0}}(f_{i}(t+x_{0},t))
e^{x_{0}\D}Y(v^{(i)},-x_{0})u^{(i)},x_{2}\right)\\
&=&\sum_{i=1}^{r}\iota_{t,x_{2},x_{0}}(f_{i}(t+x_{2}+x_{0},t+x_{2}))
Y(e^{x_{0}\D}Y(v^{(i)},-x_{0})u^{(i)},x_{2})\\
&=&\sum_{i=1}^{r}\iota_{t,x_{2},x_{0}}(f_{i}(t+x_{2}+x_{0},t+x_{2}))
Y(Y(v^{(i)},-x_{0})u^{(i)},x_{2}+x_{0}).
\end{eqnarray*}
Then
\begin{eqnarray*}
& &x_{0}^{2k}Y(Y(u,x_{0})v,x_{1}-x_{0})\\
&=&\sum_{i=1}^{r}\iota_{t,x_{1},x_{0}}(x_{0}^{k}f_{i}(t+x_{1},t+x_{1}-x_{0}))
x_{0}^{k}Y(Y(v^{(i)},-x_{0})u^{(i)},x_{1})\\
&=&\sum_{i=1}^{r}\iota_{t,x_{1},x_{0}}(x_{0}^{k}f_{i}(t+x_{1},t+x_{1}-x_{0}))\\
&&\hspace{2cm}\cdot \left((x_{1}-x_{2})^{k}Y(v^{(i)},x_{2})Y(u^{(i)},x_{1})\right)|_{x_{2}=x_{1}-x_{0}}\\
&=&\sum_{i=1}^{r}\left(\iota_{t,x_{1},x_{2}}(x_{1}-x_{2})^{k}f_{i}(t+x_{1},t+x_{2})
\right)|_{x_{2}=x_{1}-x_{0}}\\
&&\hspace{2cm}\cdot \left((x_{1}-x_{2})^{k}Y(v^{(i)},x_{2})Y(u^{(i)},x_{1})\right)|_{x_{2}=x_{1}-x_{0}}\\
&=&\left(\sum_{i=1}^{r}\iota_{t,x_{2},x_{1}}(f_{i}(t+x_{1},t+x_{2}))
(x_{1}-x_{2})^{2k}Y(v^{(i)},x_{2})Y(u^{(i)},x_{1})\right)|_{x_{2}=x_{1}-x_{0}}\\
&=&E(x_{1},x_{2})|_{x_{2}=x_{1}-x_{0}},
\end{eqnarray*}
where we are using the basic facts from Lemmas \ref{lfirst-form} and
\ref{lsubcomm}. Thus
\begin{eqnarray*}
& &\left((x_{1}-x_{2})^{2k}Y(u,x_{1})Y(v,x_{2})\right)|_{x_{1}=x_{2}+x_{0}}\\
&=&x_{0}^{2k}Y(Y(u,x_{0})v,x_{2})\\
&=&\left(E(x_{1},x_{2})|_{x_{2}=x_{1}-x_{0}}\right)|_{x_{1}=x_{2}+x_{0}}\\
&=&E(x_{1},x_{2})|_{x_{1}=x_{2}+x_{0}}.
\end{eqnarray*}
It follows that
\begin{eqnarray*}
& &(x_{1}-x_{2})^{2k}Y(u,x_{1})Y(v,x_{2})=E(x_{1},x_{2})\\
&=&\sum_{i=1}^{r}\iota_{t,x_{2},x_{1}}(f_{i}(t+x_{1},t+x_{2}))
(x_{1}-x_{2})^{2k}Y(v^{(i)},x_{2})Y(u^{(i)},x_{1}),
\end{eqnarray*}
proving (\ref{et-local}).
\end{proof}

As an immediate consequence we have:

\bc{cskew-symmetry} A nonlocal vertex $\F((t))$-algebra $V$ is a
weak quantum vertex $\F((t))$-algebra if and only if for any $u,v\in
V$, there exist
$$u^{(i)},\ v^{(i)}\in V,\; f_{i}(x_{1},x_{2})\in \F_{*}(x_{1},x_{2})
\ \ (i=1,\dots,r)$$
 such that
\begin{eqnarray*}\label{eskew-symmetry-form}
Y(u,x)v&=&\sum_{i=1}^{r}\iota_{t,x}(f_{i}(x+t,t))
e^{x\D}Y(v^{(i)},-x)u^{(i)}.
\end{eqnarray*} \ec

The following result implies that if $V$ is a weak quantum vertex
$\F((t))$-algebra, for any type zero $V$-module $W$, a variant of
$\S_{t}(x_{1},x_{2})$-Jacobi identity (\ref{et-jacobi-wqva}) holds:

\bp{pmodule-property} Let $V$ be a nonlocal vertex
$\F((t))$-algebra, let $(W,Y_{W})$ be a type zero $V$-module, and
let
$$u,\ v,\ u^{(i)},\ v^{(i)}\in V,\; f_{i}(x_{1},x_{2})\in \F_{*}(x_{1},x_{2})
\ \ (i=1,\dots,r).$$ Assume that (\ref{et-local}) holds for some
nonnegative integer $k$. Then
\begin{eqnarray}\label{et-jacobi-module}
& &x_{0}^{-1}\delta\left(\frac{x_{1}-x_{2}}{x_{0}}\right)
Y_{W}(u,x_{1})Y_{W}(v,x_{2})\nonumber\\
& &\ \ \ \ -\sum_{i=1}^{r}
x_{0}^{-1}\delta\left(\frac{x_{2}-x_{1}}{-x_{0}}\right)
\iota_{x_{2},x_{1}}(f_{i}(x_{1},x_{2}))
Y_{W}(v^{(i)},x_{2})Y_{W}(u^{(i)},x_{1})\nonumber\\
&=&x_{2}^{-1}\delta\left(\frac{x_{1}-x_{0}}{x_{2}}\right)
Y_{W}(Y(u,x_{0})v,x_{2}).
\end{eqnarray}
\ep

\begin{proof}
Since $Y_{W}(f(t)a,x)=f(x)Y_{W}(a,x)$ for $f(t)\in \F((t)),\; a\in
V$, using Proposition \ref{pskewsymmetry} and Lemma
\ref{ld-property} we get
\begin{eqnarray*}
& &Y_{W}(Y(u,x_{0})v,x_{2})\\
&=&\sum_{i=1}^{r}Y_{W}\left(\iota_{t,x_{0}}(f_{i}(t+x_{0},t))
e^{x_{0}\D}Y_{W}(v^{(i)},-x_{0})u^{(i)},x_{2}\right)\\
&=&\sum_{i=1}^{r}\iota_{x_{2},x_{0}}(f_{i}(x_{2}+x_{0},x_{2}))
Y_{W}(e^{x_{0}\D}Y(v^{(i)},-x_{0})u^{(i)},x_{2})\\
&=&\sum_{i=1}^{r}\iota_{x_{2},x_{0}}(f_{i}(x_{2}+x_{0},x_{2}))
Y_{W}(Y(v^{(i)},-x_{0})u^{(i)},x_{2}+x_{0}).
\end{eqnarray*}
Then it follows from the second half of the proof of Proposition
\ref{pskewsymmetry}.
\end{proof}

Next, we study another category of modules for nonlocal vertex
$\F((t))$-algebras.

\bd{dtmodule} {\em Let $V$ be a nonlocal vertex $\F((t))$-algebra. A
{\em type one (resp. quasi) $V$-module} is an $\F((t))$-module $W$
which is also a (resp. quasi) module for $V$ viewed as a nonlocal
vertex algebra over $\F$ such that
\begin{eqnarray}
Y_{W}(f(t)v,x)(g(t)w)=f(t+x)g(t)Y_{W}(v,x)w
\end{eqnarray}
for $f(t),g(t)\in \F((t)),\; v\in V,\; w\in W$.} \ed

We have the following simple fact:

\bl{lconsistency} Let $V$ be a nonlocal vertex $\F((t))$-algebra. a)
Let $W$ be a type one quasi $V$-module and let $U$ be a quasi
submodule of $W$ for $V$ viewed as a nonlocal vertex algebra over
$\F$. Then $U$ is an $\F((t))$-submodule of $W$. b) Let $U$ and $W$
be type one quasi $V$-modules and let $\psi: U\rightarrow W$ be a
homomorphism of quasi modules for $V$ viewed as a nonlocal vertex
algebra over $\F$. Then $\psi$ is $\F((t))$-linear. \el

\begin{proof} For a), by assumption, $U$ is
an $\F$-subspace of $W$, which is closed under the action of $V$.
For $f(t)\in \F((t)),\; w\in U$, we have
$$f(t+x)w=f(t+x)Y_{W}({\bf 1},x)w=Y_{W}(f(t){\bf 1},x)w\in U((x)),$$
which implies $f(t)w\in U$. Thus $U$ is an $\F((t))$-submodule of
$W$.

For b), we are given that $\psi$ is an $\F$-linear map such that
$$\psi(Y_{U}(v,x)u)=Y_{W}(v,x)\psi(u)\ \ \mbox{ for }v\in V,\; u\in
U.$$ For $f(t)\in \F((t)),\; u\in U$, we have
$$\psi(f(t+x)u)=\psi(Y_{U}(f(t){\bf 1},x)u)=Y_{W}(f(t){\bf
1},x)\psi(u)=f(t+x)\psi(u),$$ which implies
$\psi(f(t)u)=f(t)\psi(u)$. Thus $\psi$ is $\F((t))$-linear.
\end{proof}

The same proof (the second half) of Proposition \ref{pskewsymmetry}
yields the following analog of Proposition \ref{pmodule-property}:

\bp{pmodule-property-type2} Let $V$ be a nonlocal vertex
$\F((t))$-algebra, let $(W,Y_{W})$ be a type one $V$-module, and let
$$u,\ v,\ u^{(i)},\ v^{(i)}\in V,\; f_{i}(x_{1},x_{2})\in \F_{*}(x_{1},x_{2})
\ \ (i=1,\dots,r). $$ Assume that (\ref{et-local}) holds for some
nonnegative integer $k$. Then
\begin{eqnarray}\label{et-jacobi-module-type2}
& &x_{0}^{-1}\delta\left(\frac{x_{1}-x_{2}}{x_{0}}\right)
Y_{W}(u,x_{1})Y_{W}(v,x_{2})\nonumber\\
& &\ \ \ \ -\sum_{i=1}^{r}
x_{0}^{-1}\delta\left(\frac{x_{2}-x_{1}}{-x_{0}}\right)
\iota_{t,x_{2},x_{1}}(f_{i}(t+x_{1},t+x_{2}))
Y_{W}(v^{(i)},x_{2})Y_{W}(u^{(i)},x_{1})\nonumber\\
&=&x_{2}^{-1}\delta\left(\frac{x_{1}-x_{0}}{x_{2}}\right)
Y_{W}(Y(u,x_{0})v,x_{2}).
\end{eqnarray}
\ep

In the following two lemmas, we present some technical results which
we shall need in later sections.

\bl{lgenerating-1} Let $V$ be an $\F((t))$-module and a nonlocal
vertex algebra over $\F$. Assume that there exists a subset $U$ of
$V$ such that $\F((t))U$ generates $V$ as a nonlocal vertex algebra
over $\F$ and such that
$$Y(f(t)u,x)g(t)=f(t+x)g(t)Y(u,x)\ \ \mbox{ for }f(t),g(t)\in \F((t)),\;
u\in U\cup \{{\bf 1}\}.$$ Then $V$ is a nonlocal vertex
$\F((t))$-algebra. \el

\begin{proof} Set
$$K=\left\{ v\in V\; |\; Y(f(t)v,x)g(t)=f(t+x)g(t)Y(v,x)\ \mbox{ for }f(t),g(t)\in
\F((t))\right\}.$$ We must prove $V=K$. It is clear that $K$ is an
$\F((t))$-submodule. {}From assumption we have $\F((t))U\cup \{{\bf
1}\}\subset K$, so that $K$ generates $V$ as a nonlocal vertex
algebra. Now, it suffices to show that $K$ is closed. Let $u,v\in
K,\; f(t),g(t)\in \F((t))$. For any $w\in V$, there exists $l\in \N$
such that
\begin{eqnarray*}
&&(x_{0}+x_{2})^{l}Y(Y(f(t)u,x_{0})v,x_{2})g(t)w
=(x_{0}+x_{2})^{l} Y(f(t)u,x_{0}+x_{2})Y(v,x_{2})g(t)w\\
&&(x_{0}+x_{2})^{l}Y(u,x_{0}+x_{2})Y(v,x_{2})w
=(x_{0}+x_{2})^{l}Y(Y(u,x_{0})v,x_{2})w.
\end{eqnarray*}
Then
\begin{eqnarray*}
&&(x_{0}+x_{2})^{l}Y(f(t+x_{0})Y(u,x_{0})v,x_{2})g(t)w\\
&=&(x_{0}+x_{2})^{l}Y(Y(f(t)u,x_{0})v,x_{2})g(t)w\\
&=&(x_{0}+x_{2})^{l} Y(f(t)u,x_{0}+x_{2})Y(v,x_{2})g(t)w\\
&=&(x_{0}+x_{2})^{l}f(t+x_{0}+x_{2})Y(u,x_{0}+x_{2})Y(v,x_{2})g(t)w\\
&=&(x_{0}+x_{2})^{l}f(t+x_{0}+x_{2})g(t)Y(u,x_{0}+x_{2})Y(v,x_{2})w\\
&=& (x_{0}+x_{2})^{l}f(t+x_{0}+x_{2})g(t)Y(Y(u,x_{0})v,x_{2})w.
\end{eqnarray*}
Note that as $Y(u,x_{0})v\in V((x_{0}))$, both expressions
$$Y(f(t+x_{0})Y(u,x_{0})v,x_{2})g(t)w\ \ \mbox{and }\
f(t+x_{0}+x_{2})g(t)Y(Y(u,x_{0})v,x_{2})w$$ lie in
$V((x_{2}))((x_{0}))$. It follows that
$$Y(f(t+x_{0})Y(u,x_{0})v,x_{2})g(t)w=f(t+x_{0}+x_{2})g(t)Y(Y(u,x_{0})v,x_{2})w.$$
We note that this also holds with $f(t)$ replaced by its derivatives
of all orders.

Next, we show $u_{m}v\in K$ for $m\in \Z$ by using the above
information. Let $m\in \Z$ be arbitrarily fixed. Choosing $k\in \Z$
such that $x_{0}^{m+k}Y(u,x_{0})v\in V[[x_{0}]]$, we get
\begin{eqnarray*}
&&Y(f(t)u_{m}v,x_{2})g(t)w\\
&=&\Res_{x_{0}}x_{0}^{m}Y(f(t)Y(u,x_{0})v,x_{2})g(t)w\\
&=&\Res_{x_{0}}x_{0}^{m}Y(e^{-x_{0}\frac{\partial}{\partial t}}f(t+x_{0})Y(u,x_{0})v,x_{2})g(t)w\\
&=&\Res_{x_{0}}\sum_{n=0}^{k}\frac{(-1)^{i}}{n!}x_{0}^{m+n}Y(f^{(n)}(t+x_{0})Y(u,x_{0})v,x_{2})g(t)w\\
&=&\Res_{x_{0}}\sum_{n=0}^{k}\frac{(-1)^{i}}{n!}x_{0}^{m+n}
f^{(n)}(t+x_{0}+x_{2})g(t)Y(Y(u,x_{0})v,x_{2})w\\
&=&\Res_{x_{0}}x_{0}^{m}\left(e^{-x_{0}\frac{\partial}{\partial
t}}f(t+x_{0}+x_{2})\right)g(t)Y(Y(u,x_{0})v,x_{2})w\\
&=&f(t+x_{2})g(t)Y(u_{m}v,x_{2})w.
\end{eqnarray*}
Thus $u_{m}v\in K$. This proves that $K$ is closed, concluding the
proof.
\end{proof}

For $F(x_{1},x_{2}),G(x_{1},x_{2})\in V[[x_{1}^{\pm 1},x_{2}^{\pm
1}]]$, we define $F\sim_{\pm} G$ if
$$(x_{1}\pm x_{2})^{p}F=(x_{1}\pm x_{2})^{p}G$$
for some $p\in \N$. It is clear that the defined relations
``$\sim_{\pm}$'' are equivalence relations.

Let $U$ be a subset of a nonlocal vertex $\F((t))$-algebra $V$. We
say $U$ is {\em $\S_{t}$-local} if the $\S_{t}$-locality condition
in Proposition \ref{ptlocality} holds with $U$ in place of $V$.

We have (cf. \cite{li-qva2}, Lemma 2.7; \cite{ltw}, Proposition
2.6):

\bl{lgenerating-qva} Let $V$ be a nonlocal vertex $\F((t))$-algebra.
Assume that there exists an $\S_{t}$-local subset $U$ of $V$ such
that $\F((t))U$ generates $V$ as a nonlocal vertex algebra over
$\F$. Then $V$ is a weak quantum vertex $\F((t))$-algebra. \el

\begin{proof} First we introduce a technical notion.
 We say that an ordered
pair $(A,B)$ of subsets of $V$ is $\S_{t}$-local if for any $a\in
A,\; b\in B$, there exist
$$a^{(i)}\in A,\; b^{(i)}\in B,\; f_{i}(x_{1},x_{2})\in \F_{*}(x_{1},x_{2})
\ (i=1,\dots,r)$$ such that
$$Y(a,x_{1})Y(b,x_{2})\sim_{-} \sum_{i=1}^{r}\iota_{t,x_{2},x_{1}}(f_{i}(t+x_{1},t+x_{2}))
Y(b^{(i)},x_{2})Y(a^{(i)},x_{1}),$$
 or equivalently (in view of Corollary \ref{cskew-symmetry})
$$Y(a,x)b=\sum_{i=1}^{r}\iota_{t,x}(f_{i}(t+x,t))e^{x\D}Y(b^{(i)},-x)a^{(i)}.$$
It is clear that if $(A,B)$ is $\S_{t}$-local, so is
$(\F((t))A,\F((t))B)$. For any subset $A$ of $V$, we set
$$A^{(2)}=\F((t))\mbox{-span}\{ u_{n}v\;|\; u,v\in A,\; n\in \Z\}\subset V.$$

We are going to prove that if an ordered pair $(A,P)$ of
$\F((t))$-submodules of $V$ is $\S_{t}$-local, then $(A,P^{(2)})$
and $(A^{(2)},P)$ are $\S_{t}$-local. Then it follows {}from this
and induction that $(\<\F((t))U\>,\<\F((t))U\>)$ is $\S_{t}$-local.
Therefore, $V$ is a weak quantum vertex $\F((t))$-algebra.

We first prove that $(A,P^{(2)})$ is $\S_{t}$-local. Let $a\in A,\;
u,v\in P$. By $\S_{t}$-locality and by Proposition
\ref{pskewsymmetry}, there exist
$$f_{i}(x_{1},x_{2}), g_{ij}(x_{1},x_{2})\in \F_{*}(x_{1},x_{2}),\;
 a^{(i)}, a^{(ij)}\in A,\; u^{(i)},v^{(j)}\in P$$
 for $1\le i\le r, 1\le j\le s,$ such that
$$Y(a,x_{1})Y(u,x_{2})v\sim_{-}
\sum_{i=1}^{r}\iota_{t,x_{2},x_{1}}(f_{i}(t+x_{1},t+x_{2}))
Y(u^{(i)},x_{2})Y(a^{(i)},x_{1})v,$$
$$Y(a^{(i)},x)v=\sum_{j=1}^{s}\iota_{t,x}(g_{ij}(t+x,t))e^{x\D}Y(v^{(j)},-x)a^{(ij)},$$
 and
\begin{eqnarray}
Y(u^{(i)},x_{2}-x_{1})Y(v^{(j)},-x_{1})a^{(ij)} \sim_{-}
Y(Y(u^{(i)},x_{2})v^{(j)},-x_{1})a^{(ij)}.
\end{eqnarray}
Then using the $\D$-bracket-derivative property (\ref{ed-bracket})
and weak associativity we get
\begin{eqnarray*}\label{s-skew}
& &Y(a,x_{1})Y(u,x_{2})v\nonumber\\
&\sim_{-}
&\sum_{i=1}^{r}\iota_{t,x_{2},x_{1}}(f_{i}(t+x_{1},t+x_{2}))
Y(u^{(i)},x_{2})Y(a^{(i)},x_{1})v\nonumber\\
&\sim_{-}
&\sum_{i=1}^{r}\iota_{t,x_{2},x_{1}}(f_{i}(t+x_{1},t+x_{2}))
Y(u^{(i)},x_{2})\sum_{j=1}^{s}g_{ij}(t+x_{1},t)e^{x_{1}\D}Y(v^{(j)},-x_{1})a^{(ij)}\nonumber\\
&\sim_{-}&\sum_{i=1}^{r}\sum_{j=1}^{s}\iota_{t,x_{2},x_{1}}f_{i}(t+x_{1},t+x_{2})g_{ij}(t+x_{1},t)
e^{x_{1}\D}Y(u^{(i)},x_{2}-x_{1})Y(v^{(j)},-x_{1})a^{(ij)}\nonumber\\
&\sim_{-}&\sum_{i=1}^{r}\sum_{j=1}^{s}\iota_{t,x_{1},x_{2}}f_{i}(t+x_{1},t+x_{2})g_{ij}(t+x_{1},t)
e^{x_{1}\D}Y(Y(u^{(i)},x_{2})v^{(j)},-x_{1})a^{(ij)}.
\end{eqnarray*}
That is, there exists a nonnegative integer $k$ such that
\begin{eqnarray*}
&&(x_{1}-x_{2})^{k}Y(a,x_{1})Y(u,x_{2})v\nonumber\\
&&=(x_{1}-x_{2})^{k}\sum_{i=1}^{r}\sum_{j=1}^{s}\iota_{t,x_{1},x_{2}}
\left(f_{i}(t+x_{1},t+x_{2})g_{ij}(t+x_{1},t)\right)\\
&&\hspace{3cm}\cdot
e^{x_{1}\D}Y(Y(u^{(i)},x_{2})v^{(j)},-x_{1})a^{(ij)}.
\end{eqnarray*}
As both sides involve only finitely many negative powers of $x_{2}$,
multiplying both sides by $(x_{1}-x_{2})^{-k}$, we obtain
\begin{eqnarray*}
&&Y(a,x_{1})Y(u,x_{2})v\\
&=&\sum_{i=1}^{r}\sum_{j=1}^{s}\iota_{t,x_{1},x_{2}}f_{i}(t+x_{1},t+x_{2})g_{ij}(t+x_{1},t)
e^{x_{1}\D}Y(Y(u^{(i)},x_{2})v^{(j)},-x_{1})a^{(ij)}.
\end{eqnarray*}
It follows that $(A,P^{(2)})$ is $\S_{t}$-local.

Next, we prove that $(A^{(2)},P)$ is $\S_{t}$-local. Let $a,b\in
A,\; w\in P$. There exist $$f_{i}(x_{1},x_{2}),
g_{ij}(x_{1},x_{2})\in \F_{*}(x_{1},x_{2}),\; a^{(ij)},b^{(j)}\in
A,\; w^{(i)}, w^{(ij)}\in P$$
 for $1\le i\le
r,\; 1\le j\le s,$  such that
$$Y(b,x_{2})w=\sum_{i=1}^{r}\iota_{t,x}(f_{i}(t+x,t)) e^{x\D}Y(w^{(i)},-x)b^{(i)},$$
\begin{eqnarray*}
&&Y(a,x_{1})Y(w^{(i)},-x_{2})b^{(i)}\\
& &\sim_{+}\;
\sum_{j=1}^{s}\iota_{t,x_{2},x_{1}}(g_{ij}(t+x_{1},t-x_{2}))
Y(w^{(ij)},-x_{2})Y(a^{(ij)},x_{1})b^{(i)}.
\end{eqnarray*}
Then we get
\begin{eqnarray*}\label{s-skew2}
& &Y(Y(a,x_{1})b,x_{2})w\nonumber\\
&\sim_{+}& Y(a,x_{1}+x_{2})Y(b,x_{2})w\nonumber\\
&\sim_{+}&Y(a,x_{1}+x_{2})\sum_{i=1}^{r}\iota_{t,x_{2}}(f_{i}(t+x_{2},t))
e^{x_{2}\D}Y(w^{(i)},-x_{2})b^{(i)}\nonumber\\
&\sim_{+}&\sum_{i=1}^{r}\iota_{t,x_{2}}(f_{i}(t+x_{2},t))
e^{x_{2}\D}Y(a,x_{1})Y(w^{(i)},-x_{2})b^{(i)}\nonumber\\
&\sim_{+}&\sum_{i=1}^{r}\sum_{j=1}^{s}\iota_{t,x_{2},x_{1}}f_{i}(t+x_{2},t)
e^{x_{2}\D}g_{ij}(t+x_{1},t-x_{2})
Y(w^{(ij)},-x_{2})Y(a^{(ij)},x_{1})b^{(i)}\\
&=&\sum_{i=1}^{r}\sum_{j=1}^{s}\iota_{t,x_{2},x_{1}}(f_{i}(t+x_{2},t)g_{ij}(t+x_{1}+x_{2},t))
e^{x_{2}\D}Y(w^{(ij)},-x_{2})Y(a^{(ij)},x_{1})b^{(i)}.
\end{eqnarray*}
By a similar reasoning we obtain
\begin{eqnarray*}
&&Y(Y(a,x_{1})b,x_{2})w\\
&=&\sum_{i=1}^{r}\sum_{j=1}^{s}\iota_{t,x_{2},x_{1}}f_{i}(t+x_{2},t)
g_{ij}(t+x_{1}+x_{2},t)e^{x_{2}\D}
Y(w^{(ij)},-x_{2})Y(a^{(ij)},x_{1})b^{(i)}.
\end{eqnarray*}
It follows that $(A^{(2)},P)$ is $\S_{t}$-local. Now, the proof is
complete.
\end{proof}

\section{Quantum vertex $\F((t))$-algebras and non-degeneracy}
In this section we formulate and study a notion of quantum vertex
$\F((t))$-algebra and we study Etingof-Kazhdan's notion of
non-degeneracy for nonlocal vertex $\F((t))$-algebras. As a key
result we show that every non-degenerate weak quantum vertex
$\F((t))$-algebra has a canonical quantum vertex $\F((t))$-algebra
structure. In this section we also present some basic results on
non-degeneracy.

We begin with some basics on quantum Yang-Baxter operators. Let $H$
be a vector space over $\F$. The symmetric group $S_{3}$ naturally
acts on $H^{\otimes 3}$ with $\sigma\in S_{3}$ acting as
$P_{\sigma}$ which is defined by
$$P_{\sigma}(u_{1}\otimes u_{2}\otimes u_{3})=u_{\sigma(1)}\otimes
u_{\sigma(2)}\otimes u_{\sigma(3)}\ \ \ \mbox{ for  }
u_{1},u_{2},u_{3}\in H.$$ For $1\le i<j\le 3$, set $P_{ij}=P_{(ij)}$
(with $(ij)$ denoting the transposition). We have
$$P_{12}P_{23}P_{12}=P_{13}=P_{23}P_{12}P_{23}.$$
Let $P$ denote the flip operator on $H\otimes H$ with $P(u\otimes
v)=v\otimes u$ for $u,v\in H$. Then
$$P_{12}=P\otimes 1,\ \ \ \ P_{23}=1\otimes P.$$
A {\em quantum Yang-Baxter operator with two parameters} on $H$ is a
linear map
$$\S(x_{1},x_{2}): H\otimes H\rightarrow H\otimes H\otimes
\F_{*}(x_{1},x_{2}),$$ satisfying the quantum Yang-Baxter equation
\begin{eqnarray}\label{egqyb}
\S_{12}(x_{1},x_{2})\S_{13}(x_{1},x_{3})\S_{23}(x_{2},x_{3})=
\S_{23}(x_{2},x_{3})\S_{13}(x_{1},x_{3})\S_{12}(x_{1},x_{2}),
\end{eqnarray}
where $\S_{ij}(x_{i},x_{j})$ are the linear maps from $H^{\otimes
3}\rightarrow H^{\otimes 3}\otimes \F_{*}(x_{i},x_{j})$, defined by
$\S_{12}(x,z)=\S(x,z)\otimes 1,\ \ \S_{23}(x,z)=1\otimes \S(x,z)$,
and
$$\S_{13}(x,z)=P_{23}(\S(x,z)\otimes 1)P_{23}.$$
Furthermore, $\S(x_{1},x_{2})$ is said to be {\em unitary} if
\begin{eqnarray}
\S_{21}(x_{2},x_{1})\S(x_{1},x_{2})=1,
\end{eqnarray}
where $\S_{21}(x_{2},x_{1})=P\S(x_{2},x_{1})P$. Set
\begin{eqnarray}
R(x_{1},x_{2})=\S(x_{1},x_{2})P:\ \ H\otimes H\rightarrow H\otimes
H\otimes \F_{*}(x_{1},x_{2}).
\end{eqnarray}
It is known that (\ref{egqyb}) is equivalent to the following
braided relation
\begin{eqnarray}
R_{12}(x_{1},x_{2})R_{23}(x_{1},x_{3})R_{12}(x_{2},x_{3})
=R_{23}(x_{2},x_{3})R_{12}(x_{1},x_{3})R_{23}(x_{1},x_{2}).
\end{eqnarray}

\bd{dtqva-s} {\em A {\em quantum vertex $\F((t))$-algebra} is a weak
quantum vertex $\F((t))$-algebra $V$ equipped with an $\F$-linear
unitary quantum Yang-Baxter operator $\S(x_{1},x_{2})$ (with two
parameters) on $V$, satisfying the conditions that
\begin{eqnarray}\label{esemi-linear}
\S(x_{1},x_{2})(f(t)u\otimes
g(t)v)=f(x_{1})g(x_{2})\S(x_{1},x_{2})(u\otimes v)
\end{eqnarray}
for $f(t),g(t)\in \F((t)),\; u,v\in V$, and that for $u,v\in V$,
\begin{eqnarray*}
& &(x_{1}-x_{2})^{k}Y(v,x_{2})Y(u,x_{1})\nonumber\\
&= &\sum_{i=1}^{r}(x_{1}-x_{2})^{k}
\iota_{t,x_{1},x_{2}}(f_{i}(x_{1}+t,x_{2}+t))
Y(u^{(i)},x_{1})Y(v^{(i)},x_{2})
\end{eqnarray*}
for some nonnegative integer $k$, where $u^{(i)},v^{(i)}, f_{i}\
(i=1,\dots,r)$ are given by
$$\S(x_{1},x_{2})(u\otimes v)
=\sum_{i=1}^{r}u^{(i)}\otimes v^{(i)}\otimes f_{i}(x_{1},x_{2}), $$
and that
\begin{eqnarray}
& &[\D\otimes 1,\S(x_{1},x_{2})]=-\frac{\partial}{\partial
x_{1}}\S(x_{1},x_{2}), \ \  [1\otimes \D, \S(x_{1},x_{2})]
=-\frac{\partial}{\partial x_{2}}\S(x_{1},x_{2}),\ \ \ \ \ \ \ \
\label{edbracket-der}\\
& &\S(x_{1},x_{2})(Y(x)\otimes 1) =(Y(x)\otimes
1)\S_{23}(x_{1},x_{2})\S_{13}(x_{1}+x,x_{2}).\label{ehexagon}
\end{eqnarray}}
\ed

We modify Etingof-Kazhdan's notion of non-degeneracy (see \cite{ek})
as follows:

\bd{dnondegenerate-ta} {\em Let $V$ be a nonlocal vertex
$\F((t))$-algebra. Denote by $V^{\otimes n}$ the tensor product
space over $\F$ and define $ V^{\otimes n}\boxtimes
\F_{*}(x_{1},\dots,x_{n})$ to be the quotient space of $V^{\otimes
n}\otimes \F_{*}(x_{1},\dots,x_{n})$ by the relations
$$ f_{1}(t)v^{(1)}\otimes \cdots \otimes f_{n}(t)v^{(n)}\otimes f
=  v^{(1)}\otimes \cdots \otimes v^{(n)}\otimes f_{1}(x_{1})\cdots
f_{n}(x_{n})f
$$
for $f\in \F_{*}(x_{1},\dots,x_{n}),\; f_{i}(t)\in \F((t)),\;
v^{(i)}\in V$ $(i=1,\dots,n)$. We say that $V$ is {\em
non-degenerate} if for every positive integer $n$, the $\F$-linear
map
\begin{eqnarray*}
Z_{n}: V^{\otimes n}\boxtimes  \F_{*}(x_{1},\dots,x_{n}) &
&\rightarrow V((x_{1}))\cdots ((x_{n}))\\
(v^{(1)}\otimes \cdots \otimes v^{(n)})\boxtimes f & & \mapsto
\iota_{t,x_{1},\dots,x_{n}}f(t+x_{1},\dots, t+x_{n})
Y(v^{(1)},x_{1})\cdots Y(v^{(n)},x_{n}){\bf 1}
\end{eqnarray*}
is injective. (One can see that $Z_{n}$ is indeed well defined.)}
\ed

\br{rnon-degeneracy} {\em Given a nonlocal vertex $\F((t))$-algebra
$V$, let $V^{0}$ be an $\F$-subspace such that
$V=\F((t))\otimes_{\F}V^{0}$. We see that $Z_{n}$ is injective if
and only if the restriction
$$Z_{n}^{0}:\ \ (V^{0})^{\otimes n}\otimes\F_{*}(x_{1},\dots,
x_{n})  \rightarrow V((x_{1}))\cdots ((x_{n}))$$ is injective. For
$g_{i}(x)\in \F((x)),\; v^{(i)}\in V^{0}$ $(i=1,\dots,r)$, we have
\begin{eqnarray*}
&&Z_{1}\left(\sum_{i=1}^{r}v^{(i)}\otimes g_{i}(x)
\right)=\sum_{i=1}^{r}g_{i}(t+x)Y(v^{(i)},x){\bf
1}=\sum_{i=1}^{r}g_{i}(t+x)e^{x\D}v^{(i)}\\
&&\hspace{3cm}=e^{x\D}\sum_{i=1}^{r}g_{i}(t)v^{(i)}.
\end{eqnarray*}
{}From this we see that $Z_{1}$ is always injective.} \er

Actually, what we need in practice are certain variations of the
maps $Z_{n}$.

\bd{dpinmaps} {\em For each $n\ge 1$, we define an $\F$-linear map
\begin{eqnarray}
\pi_{n}:\ \ V^{\otimes n}\boxtimes\F_{*}(x_{1},\dots,x_{n})
\rightarrow \Hom (V,V((x_{1}))\cdots ((x_{n})))
\end{eqnarray}
by
\begin{eqnarray*}
&&\pi_{n}(v^{(1)}\otimes \cdots \otimes
v^{(n)}\otimes f)\\
&=&\iota_{t,x_{1},\dots,x_{n}}f(t+x_{1},\dots, t+x_{n})
Y(v^{(1)},x_{1})\cdots Y(v^{(n)},x_{n})
\end{eqnarray*}
for $f\in \F_{*}(x_{1},\dots,x_{n}),\; v^{(1)},\dots,v^{(n)}\in V$.}
\ed

Noticing that
$$\pi_{n}(v^{(1)}\otimes \cdots \otimes
v^{(n)}\otimes f)({\bf 1})=Z_{n}(v^{(1)}\otimes \cdots \otimes
v^{(n)}\otimes f),$$  we see that the injectivity of $Z_{n}$ implies
the injectivity of $\pi_{n}$.

We follow \cite{ek} to denote by
$$Y(x):\ \ V\otimes V\rightarrow V((x))\subset V[[x,x^{-1}]]$$
the $\F$-linear map defined by $Y(x)(u\otimes v)=Y(u,x)v$ for
$u,v\in V$. As a common practice, $Y(x)$ is always extended:
$$Y(x): V\otimes V\otimes
\F_{*}(t,x_{1},\dots,x_{k})((x))\rightarrow (V\otimes
\F_{*}(t,x_{1},\dots,x_{k}))((x)),$$ where
\begin{eqnarray}
Y(x)(u\otimes v\otimes f)=fY(x)(u\otimes v)=fY(u,x)v
\end{eqnarray}
for $u,v\in V,\; f\in \F_{*}(t,x_{1},\dots,x_{k})((x))$, where $k$
is a positive integer.

The following, which is lifted from \cite{ek} (Proposition 1.11),
plays a very important role in the theory of quantum vertex
$\F((t))$-algebras:

\bt{tek} Let $V$ be a weak quantum vertex $\F((t))$-algebra. Assume
that $V$ is non-degenerate. Then there exists an $\F$-linear map
$$\S(x_{1},x_{2}): V\otimes V\rightarrow V\otimes V\otimes
\F_{*}(x_{1},x_{2}),$$ which is uniquely determined by the condition
that for $u,v\in V$,
\begin{eqnarray}
Y(v,x_{2})Y(u,x_{1})w\sim_{-} Y(x_{1})(1\otimes
Y(x_{2}))(\S(t+x_{1},t+x_{2})(u\otimes v)\otimes w)
\end{eqnarray}
for all $w\in V$.  Furthermore, $\S(x_{1},x_{2})$ is a unitary
quantum Yang-Baxter operator on $V$, and $V$ equipped with
$\S(x_{1},x_{2})$ is a quantum vertex $\F((t))$-algebra. \et

\begin{proof} First of all, with $V$ non-degenerate, all the maps
$\pi_{n}$ $(n\ge 1)$ are injective. Notice that
$$Y(x_{1})(1\otimes
Y(x_{2}))(\S(t+x_{1},t+x_{2})(u\otimes v)\otimes
w)=\left(\pi_{2}\S(x_{1},x_{2})(u\otimes v)\right)(w).$$ It follows
that $\S(x_{1},x_{2})$ is uniquely determined by the very condition.

Let $u,v,w\in V,\; f(t),g(t)\in \F((t))$. We have
\begin{eqnarray*}
&&Y(x_{1})(1\otimes Y(x_{2}))\S_{12}(t+x_{1},t+x_{2}) (f(t)u\otimes
g(t)v\otimes w),\\
& \sim_{-}&Y(x_{2})(1\otimes Y(x_{1}))(g(t) v\otimes f(t)u\otimes w)\\
&=&f(t+x_{1})g(t+x_{2})Y(x_{2})(1\otimes Y(x_{1}))(v\otimes u\otimes
w)\\
&\sim_{-}&f(t+x_{1})g(t+x_{2})Y(x_{1})(1\otimes
Y(x_{2}))\S_{12}(t+x_{1},t+x_{2})(u\otimes v\otimes w)\\
&=&Y(x_{1})(1\otimes
Y(x_{2}))f(t+x_{1})g(t+x_{2})\S_{12}(t+x_{1},t+x_{2})(u\otimes
v\otimes w).
\end{eqnarray*}
As the first term and the last term both lie in
$V((x_{1}))((x_{2}))$, the equivalence relation between them
actually amounts to equality. With $\pi_{2}$ injective, we obtain
$$\S(x_{1},x_{2})(f(t)u\otimes
g(t)v)=f(x_{1})g(x_{2})\S(x_{1},x_{2})(u\otimes v).$$

The quantum Yang-Baxter relation, the unitarity, and the
$\D$-bracket-derivative property (\ref{edbracket-der}) follow from
the same proof of Theorem 4.8 in \cite{li-qva1} with obvious
modifications. It remains to prove (\ref{ehexagon}). For $u,v,w,a\in
V$, we have
\begin{eqnarray*}
&&\Res_{x}x^{n}Y(x_{1})(1\otimes
Y(x_{2}))\left(\S(x_{1}+t,x_{2}+t)(Y(u,x)v\otimes w)\otimes
a\right)\\
&\sim_{-}&\Res_{x}x^{n}Y(w,x_{2})Y(Y(u,x)v,x_{1})a
\end{eqnarray*}
for any fixed $n\in \Z$. On the other hand, we have
\begin{eqnarray}\label{eotoh}
&&Y(w,x_{2})Y(u,z)Y(v,x_{1})a\nonumber\\
&\sim_{-}&Y(z)(1\otimes Y(x_{2}))(1\otimes 1\otimes
Y(x_{1}))\S_{12}(z+t,x_{2}+t)(u\otimes w\otimes v\otimes a)\nonumber\\
&\sim_{-}& Y(z)(1\otimes Y(x_{1}))(1\otimes 1\otimes
Y(x_{2}))\S_{23}(x_{1}+t,x_{2}+t)P_{23}\cdot\nonumber\\
&&\ \ \ \ \cdot \S_{12}(z+t,x_{2}+t)P_{23}(u\otimes
v\otimes w\otimes a)\nonumber\\
&=& Y(z)(1\otimes Y(x_{1}))(1\otimes 1\otimes
Y(x_{2}))\cdot\nonumber\\
&&\ \ \ \ \cdot
\S_{23}(x_{1}+t,x_{2}+t)\S_{13}(z+t,x_{2}+t)(u\otimes v\otimes
w\otimes a).
\end{eqnarray}
Notice that for any $u',v'\in V$, there exists a
nonnegative integer $k$ such that
$$(x_{1}-x_{2})^{k}Y(u',x_{1})Y(v',x_{2})\in
\Hom(V,V((x_{1},x_{2}))),$$
$$x_{0}^{k}Y(Y(u',x_{0})v',x_{2})=\left(
(x_{1}-x_{2})^{k}Y(u',x_{1})Y(v',x_{2})\right)|_{x_{1}=x_{2}+x_{0}}.$$
Using (\ref{eotoh}), by choosing $k$ sufficiently large, we have
\begin{eqnarray*}
&&x^{k}Y(w,x_{2})Y(Y(u,x)v,x_{1})a\\
&=&\left((z-x_{1})^{k}Y(w,x_{2})Y(u,z)Y(v,x_{1})a\right)|_{z=x_{1}+x}\\
&\sim_{-}&x^{k}Y(x_{1})(Y(x)\otimes 1)(1\otimes 1\otimes
Y(x_{2}))\cdot\\
&&\ \ \ \ \cdot
\S_{23}(x_{1}+t,x_{2}+t)\S_{13}(x_{1}+x+t,x_{2}+t)(u\otimes v\otimes
w\otimes a).
\end{eqnarray*}
Thus
\begin{eqnarray*}
&& \Res_{x}x^{n}Y(x_{1})(1\otimes
Y(x_{2}))\left(\S(x_{1}+t,x_{2}+t)(Y(u,x)v\otimes w)\otimes
a\right)\\ &\sim_{-}&\Res_{x}x^{n}Y(x_{1})(Y(x)\otimes 1)(1\otimes
1\otimes Y(x_{2}))\cdot\\
&&\ \ \ \ \cdot
\S_{23}(x_{1}+t,x_{2}+t)\S_{13}(x_{1}+x+t,x_{2}+t)(u\otimes v\otimes
w\otimes a)
\end{eqnarray*}
for any fixed $n\in \Z$. As both sides are in
$V((x_{1}))((x_{2},x))$, we have
\begin{eqnarray*}
&& \Res_{x}x^{n}Y(x_{1})(1\otimes
Y(x_{2}))\left(\S(x_{1}+t,x_{2}+t)(Y(u,x)v\otimes w)\otimes
a\right)\\ &=&\Res_{x}x^{n}Y(x_{1})(Y(x)\otimes 1)(1\otimes
1\otimes Y(x_{2}))\cdot\\
&&\ \ \ \ \cdot
\S_{23}(x_{1}+t,x_{2}+t)\S_{13}(x_{1}+x+t,x_{2}+t)(u\otimes v\otimes
w\otimes a)\\
&=&\Res_{x}x^{n}Y(x_{1})(1\otimes Y(x_{2}))(Y(x)\otimes 1\otimes 1)\cdot\\
&&\ \ \ \ \cdot
\S_{23}(x_{1}+t,x_{2}+t)\S_{13}(x_{1}+x+t,x_{2}+t)(u\otimes v\otimes
w\otimes a).
\end{eqnarray*}
Since $n$ is arbitrary, we can drop off $\Res_{x}x^{n}$. Then
(\ref{ehexagon}) follows.
\end{proof}

For the rest of this section we focus on non-degeneracy of nonlocal
vertex $\F((t))$-algebras. Let $V$ be a nonlocal vertex
$\F((t))$-algebra. {}From Lemma \ref{lconsistency}, a $V$-submodule
of $V$ for $V$ viewed as a nonlocal vertex algebra over $\F$ is the
same as a $V$-submodule of $V$ for $V$ viewed as a nonlocal vertex
$\F((t))$-algebra.  Furthermore, a module endomorphism for $V$
viewed as a nonlocal vertex algebra over $\F$ is the same as a
module endomorphism for $V$ viewed as a nonlocal vertex
$\F((t))$-algebra. We denote by $V^{\mod}$ the adjoint $V$-module.

\bp{pnon-degenerate} Let $V$ be a nonlocal vertex $\F((t))$-algebra
such that $V$ as a $V$-module is irreducible with
$\End_{V}(V^{\mod})=\F((t))$. Then $V$ is non-degenerate. \ep

\begin{proof} We are going to use induction to show that $Z_{n}$ is injective for
every positive integer $n$, following the proof of a similar result
in \cite{li-qva2}. Recall from Remark \ref{rnon-degeneracy} that
$Z_{1}$ is always injective. Now, assume that $n\ge 2$ and $Z_{n-1}$
is injective. Let $U$ be the quotient space of $V^{\otimes
(n-1)}\otimes \F_{*}(x_{1},\dots,x_{n})$, viewed as a vector space
over $\F$, by the relations
\begin{eqnarray*}
&&f_{2}(t)v^{(2)}\otimes \cdots \otimes
f_{n}(t)v^{(n)}\otimes f\\
&=&  v^{(2)}\otimes \cdots \otimes v^{(n)}\otimes f_{2}(x_{2})\cdots
f_{n}(x_{n})f
\end{eqnarray*}
for $f\in \F_{*}(x_{1},\dots,x_{n}),\; f_{i}(t)\in \F((t)),\;
v^{(i)}\in V$ $(i=2,\dots,n)$. Note that $U$ is naturally an
$\F_{*}(x_{1},\dots,x_{n})$-module while $\F_{*}(x_{1},\dots,x_{n})$
is an algebra over $\F((x_{1}))$. Furthermore, viewing $V$ as an
$\F((x_{1}))$-module with $f(x_{1})$ acting as $f(t)$, we have
$$V^{\otimes n}\boxtimes\F_{*}(x_{1},\dots,x_{n})
=V\otimes_{\F((x_{1}))} U.$$ Let $B$ be the subalgebra of the
endomorphism algebra $\End_{\F((t))}(V)$ (over $\F((t))$) generated
by $v_{n}$ for $v\in V,\; n\in \Z$. Then $V$ is an irreducible
$B$-module with $\End_{B}(V)=\F((t))$. {}From \cite{li-qva2} (Lemma
3.8), the kernel of $Z_{n}$ is a $B\otimes_{\F((t_{1}))}
\F_{*}(x_{1},\dots,x_{n})$-submodule of $V\otimes_{\F((x_{1}))} U$
(with $B$ acting on the first factor). By a classical fact (cf.
\cite{li-simple}, Lemma 2.10), we have $\ker
Z_{n}=V\otimes_{\F((x_{1}))}P$ for some submodule $P$ of $U$. Let
$a\in P\subset U$. There exists a nonzero polynomial
$q(x_{1},\dots,x_{n})$ such that
 $$q(x_{1},\dots,x_{n})a\in V^{\otimes (n-1)}\otimes
 \F[[x_{1},\dots,x_{n}]].$$
Write
$$q(x_{1},\dots,x_{n})a=\sum_{m\in \N}x_{1}^{m}a_{m}$$
with $a_{m}\in V^{\otimes (n-1)}\otimes \F[[x_{2},\dots,x_{n}]].$ As
${\bf 1}\otimes qa\in \ker Z_{n}$, we have $a_{m}\in \ker Z_{n-1}$
for $m\in \Z$. Then $a_{m}=0$ for $m\in \Z$, and hence
$q(x_{1},\dots,x_{n})a=0$. Thus $a=0$. This proves that $P=0$, which
implies that $Z_{n}$ is injective.
\end{proof}

\br{rfiltration-graded} {\em Let $V$ be a nonlocal vertex
$\F((t))$-algebra and let ${\mathcal{F}}=\{F_{n}\}_{n\in
\frac{1}{2}\Z}$ be an increasing filtration of $\F((t))$-submodules
of $V$, satisfying the condition that ${\bf 1}\in F_{0}$,
$$u_{k}F_{n}\subset F_{m+n-k-1}\ \ \mbox{ for }u\in
F_{m},\; k\in \Z,\;  m,n\in \frac{1}{2}\Z.$$ Form the
$\frac{1}{\Z}$-graded $\F((t))$-module
$${\rm Gr}_{\mathcal{F}}(V)=\oplus_{n\in \frac{1}{2}\Z}(F_{n}/F_{n-1/2}).$$
For $u+F_{m-1/2}\in (F_{m}/F_{m-1/2}),\; v+F_{n-1/2}\in
(F_{n}/F_{n-1/2})$ with $m,n\in \frac{1}{2}\Z$, define
$$(u+F_{m-1/2})_{k}(v+F_{n-1/2})=u_{k}v+F_{m+n-k-3/2}\in (F_{m+n-k-1}/F_{m+n-k-3/2})$$
for $k\in \Z$. It is straightforward to show that ${\rm
Gr}_{\mathcal{F}}(V)$ is a nonlocal vertex $\F((t))$-algebra with
${\bf 1}+F_{-1/2}$ as the vacuum vector (cf. \cite{kl}).}\er

\bp{pfiltration-f} Let $V$ be a nonlocal vertex $\F((t))$-algebra
with an increasing filtration ${\mathcal{F}}=\{F_{n}\}_{n\in
\frac{1}{2}\Z}$ of $\F((t))$-submodules, satisfying the condition
that $F_{n}=0$ for $n$ sufficiently negative, ${\bf 1}\in F_{0}$,
and
$$u_{k}F_{n}\subset F_{m+n-k-1}\ \ \mbox{ for }u\in
F_{m}, \; k\in \Z,\; m,n\in \frac{1}{2}\Z.$$ Assume that ${\rm
Gr}_{\mathcal{F}}(V)$ as a ${\rm Gr}_{\mathcal{F}}(V)$-module is
irreducible with $\End ({\rm Gr}_{\mathcal{F}}(V)^{\mod})=\F((t))$.
Then $V$ as a $V$-module is irreducible with
$\End_{V}(V^{\mod})=\F((t))$ and $V$ is non-degenerate.
 \ep

\begin{proof} Notice that the assertion on non-degeneracy follows
{}from the other assertions and Proposition \ref{pnon-degenerate}.
The irreducibility assertion follows from Proposition 2.11 of
\cite{kl}. It remains to prove $\End_{V}(V^{\mod})=\F((t))$. Let
$\psi\in \End_{V}(V^{\mod})$. If $\psi({\bf 1})=0$, we have
$\psi=0\in \F((t))$ as
$$\psi(v)=\psi(v_{-1}{\bf 1})=v_{-1}\psi({\bf 1})=0\ \ \mbox{ for
}v\in V.$$ Thus, $\psi({\bf 1})\ne 0$ for any nonzero $\psi\in
\End_{V}(V^{\mod})$. Assume $\psi\ne 0$. Since $\psi({\bf 1})\ne 0$
and since $F_{n}=0$ for $n$ sufficiently negative, there exists
$m\in \frac{1}{2}\Z$ such that $\psi({\bf 1})\in F_{m}-F_{m-1/2}$.
For $v\in V$, we have
$$v_{n}\psi({\bf 1})=\psi(v_{n}{\bf 1})=0\ \ \mbox{ for }
n\ge 0.$$ By Lemma 6.1 of \cite{li-qva1}, we have a ${\rm
Gr}_{\mathcal{F}}(V)$-module endomorphism $\bar{\psi}$ of ${\rm
Gr}_{\mathcal{F}}(V)$, sending ${\bf 1}+F_{-1/2}\in F_{0}/F_{-1/2}$
to $\psi({\bf 1})+F_{m}/F_{m-1/2}$. {}From assumption we have
$\bar{\psi}=f(t)$ for some $f(t)\in \F((t))$. As ${\rm
Gr}_{\mathcal{F}}(V)$ is $\frac{1}{2}\Z$-graded, we must have that
$m=0$ and $\psi({\bf 1})-f(t){\bf 1}\in F_{-1/2}$. If $\psi\ne
f(t)$, with $\psi-f(t)$ in place of $\psi$ we have $(\psi-f(t))({\bf
1})\in F_{0}-F_{-1/2}$, a contradiction. Thus, $\psi=f(t)\in
\F((t))$.
\end{proof}

\br{rfiltraionE} {\em Let $V$ be a nonlocal vertex $\F((t))$-algebra
and let $\E=\{ E_{n}\}_{n\in \Z}$ be an increasing filtration of
$\F((t))$-submodules of $V$, satisfying the condition that ${\bf
1}\in E_{0}$,
$$u_{k}E_{n}\subset E_{m+n}\ \ \ \mbox{ for }u\in
E_{m},\; m,n,k\in \Z.$$ Form the $\Z$-graded $\F((t))$-module
 $${\rm Gr}_{\E}(V)=\oplus_{n\in \Z}(E_{n}/E_{n-1}).$$
For $u\in E_{m},\; v\in E_{n}$ with $m,n\in \Z$ and for $k\in \Z$,
define
$$(u+E_{m-1})_{k}(v+E_{n-1})=u_{k}v+E_{m+n-1}\in
(E_{m+n}/E_{m+n-1}).$$ It is straightforward to show that ${\rm
Gr}_{\E}(V)$ is a nonlocal vertex $\F((t))$-algebra with ${\bf
1}+E_{-1}\in E_{0}/E_{-1}$ as the vacuum vector (cf.
\cite{li-qva2}).} \er

The following follows from the same proof of Proposition 3.14 of
\cite{li-qva2} (with obvious notational modifications):

\bp{pfiltration-nondeg} Let $V$ be a nonlocal vertex
$\F((t))$-algebra and let $\E=\{ E_{n}\}_{n\in \Z}$ be an increasing
filtration of $\F((t))$-submodules, satisfying the condition that
$E_{n}=0$ for $n$ sufficiently negative, ${\bf 1}\in E_{0}$, and
$$u_{k}E_{n}\subset E_{m+n}\ \ \ \mbox{ for } u\in
E_{m},\; m,n,k\in \Z.$$ If ${\rm Gr}_{\E}(V)$ is non-degenerate,
then $V$ is non-degenerate. \ep

Let $U$ be a nonlocal vertex $\F((t))$-algebra and let $K$ be a
nonlocal vertex algebra over $\F$. Equip $U\otimes K$ with the
$\F((t))$-module structure with $\F((t))$ acting on $U$ and also
equip $U\otimes K$ with the nonlocal vertex algebra structure by
tensor product over $\F$. It can be readily seen that $U\otimes K$
becomes a nonlocal vertex $\F((t))$-algebra. Note that from
Borcherds' construction of vertex algebras, $\F((t))$ is a vertex
algebra with $1$ as the vacuum vector and with
$$Y(f(t),x)g(t)=(e^{x\frac{d}{dt}}f(t))g(t)=f(t+x)g(t)$$
for $f(t),g(t)\in \F((t))$. Thus, for any nonlocal vertex algebra
$V^{0}$ over $\F$, $\F((t))\otimes V^{0}$ is a nonlocal vertex
$\F((t))$-algebra. We have:

\bl{lnondegeneracy} Let $V^{0}$ be a non-degenerate nonlocal vertex
algebra over $\F$. Then the nonlocal vertex $\F((t))$-algebra
$\F((t))\otimes V^{0}$ is non-degenerate.
 \el

\begin{proof} {}From Remark \ref{rnon-degeneracy}, for $n\ge 1$, $Z_{n}$ is
injective if and only if the restriction
$$Z_{n}^{0}:
(V^{0})^{\otimes n}\otimes \F_{*}(x_{1},\dots,x_{n})\rightarrow
(\F((t))\otimes V^{0})((x_{1}))\cdots ((x_{n}))$$ is injective.
Furthermore, we see that $Z_{n}^{0}$ is injective if and only if its
restriction on $(V^{0})^{\otimes n}\otimes \F[[x_{1},\dots,x_{n}]]$
is injective. Assume that $A\in (V^{0})^{\otimes n}\otimes
\F[[x_{1},\dots,x_{n}]]$ such that $Z_{n}^{0}(A)=0$. By extracting
the constant term in variable $t$, we see that $A$ lies in the
kernel of the $Z_{n}$-map for $V^{0}$. Then it follows.
\end{proof}

\section{Conceptual construction of weak quantum vertex
$\F((t))$-algebras and their modules}

In this section, we present a conceptual construction of nonlocal
vertex $\F((t))$-algebras, weak quantum vertex $\F((t))$-algebras,
and their quasi modules of type zero, by using quasi compatible
subsets and quasi $\S(x_{1},x_{2})$-local subsets of formal vertex
operators. This construction is based on the conceptual construction
in \cite{li-qva1} of nonlocal vertex algebras and their quasi
modules.

We begin with the conceptual construction of nonlocal vertex
algebras and their (quasi) modules, established in \cite{li-qva1}.
Let $W$ be a vector space over $\F$. Set
\begin{eqnarray}
\E(W)=\Hom (W,W((x)))\ \ (\subset (\End W)[[x,x^{-1}]]),
\end{eqnarray}
which contains the identity operator $1_{W}$ on $W$ as a special
element.

\bd{dqcompatible} {\em A finite sequence $a_{1}(x),\dots,a_{r}(x)$
in $\E(W)$ is said to be {\em quasi compatible} if there exists a
nonzero polynomial $p(x,y)\in \F[x,y]$ such that
\begin{eqnarray}
\left(\prod_{1\le i<j\le r}p(x_{i},x_{j})\right) a_{1}(x_{1})\cdots
a_{r}(x_{r})\in \Hom (W,W((x_{1},\dots,x_{r}))).
\end{eqnarray}
The sequence $a_{1}(x),\dots,a_{r}(x)$
 is said to be {\em compatible} if there exists a nonnegative integer
$k$ such that
\begin{eqnarray}
\left(\prod_{1\le i<j\le r}(x_{i}-x_{j})^{k}\right)
a_{1}(x_{1})\cdots a_{r}(x_{r})\in \Hom (W,W((x_{1},\dots,x_{r}))).
\end{eqnarray}
Furthermore, a subset $T$ of $\E(W)$ is said to be {\em quasi
compatible} (resp. {\em compatible})  if every finite sequence in
$T$ is quasi compatible (resp. compatible).} \ed

Let $(a(x),b(x))$ be a quasi compatible ordered pair in $\E(W)$. By
definition, there exists a nonzero polynomial $p(x,y)\in \F[x,y]$
such that
\begin{eqnarray}\label{ehalf-relation}
p(x_{1},x_{2})a(x_{1})b(x_{2})\in \Hom (W,W((x_{1},x_{2}))).
\end{eqnarray}
Define $a(x)_{n}b(x)\in \E(W)$ for $n\in \Z$ in terms of the
generating function
\begin{eqnarray}
Y_{\E}(a(x),x_{0})b(x)=\sum_{n\in \Z}a(x)_{n}b(x) x_{0}^{-n-1}
\end{eqnarray}
 by
\begin{eqnarray}
Y_{\E}(a(x),x_{0})b(x)=\iota_{x,x_{0}}\left(\frac{1}{p(x+x_{0},x)}\right)
\left(p(x_{1},x)a(x_{1})b(x)\right)|_{x_{1}=x+x_{0}}.
\end{eqnarray}
A quasi compatible $\F$-subspace $U$ of $\E(W)$ is said to be {\em
$Y_{\E}$-closed} if
\begin{eqnarray*}
a(x)_{n}b(x)\in U\ \ \ \mbox{ for }a(x),b(x)\in U,\; n\in \Z.
\end{eqnarray*}

The following was obtained in \cite{li-qva1} (though the scalar
field therein is $\C$, it is clear that the results hold for any
field of characteristic $0$):

\bt{tquasi-main} Let $W$ be a vector space over $\F$ and let $U$ be
a (resp. quasi) compatible subset of $\E(W)$. There exists a
$Y_{\E}$-closed (resp. quasi) compatible subspace that contains $U$
and $1_{W}$. Denote by $\<U\>$ the smallest such subspace of
$\E(W)$. Then $(\<U\>,Y_{\E},1_{W})$ carries the structure of a
nonlocal vertex algebra with $W$ as a (resp. quasi) module where
$Y_{W}(\alpha(x),x_{0})=\alpha(x_{0})$ for $\alpha(x)\in \<U\>$. \et

Let $W$ be a vector space over $\F$ as before. Notice that for any
$f(x)\in \F((x)),\; a(x)\in \E(W)\ (=\Hom (W,W((x))))$, we have
$f(x)a(x)\in \E(W)$. Thus, $\E(W)$ is naturally an $\F((x))$-module,
namely, a vector space over the field $\F((x))$.

We now present our first main result of this section.

\bt{tconcrete-general} Let $W$ be a vector space over $\F$ and let
$U$ be any (resp. quasi) compatible subset of $\E(W)$. Let $\<U\>$
be the $Y_{\E}$-closed (resp. quasi) compatible $\F$-subspace of
$\E(W)$ as in Theorem \ref{tquasi-main}. Then $\F((x))\<U\>$ is a
$Y_{\E}$-closed (resp. quasi) compatible $\F((x))$-submodule of
$\E(W)$. Furthermore, $(\F((x))\<U\>,Y_{\E},1_{W})$ carries the
structure of a nonlocal vertex $\F((t))$-algebra, where
$$f(t)a(x)=f(x)a(x)\ \ \ \mbox{ for }f(t)\in \F((t)),\; a(x)\in
\F((x))\<U\>,$$ and $(W,Y_{W})$ carries the structure of a (resp.
quasi) $\F((x))\<U\>$-module of type zero with
$Y_{W}(a(x),x_{0})=a(x_{0})$ for $a(x)\in \F((x))\<U\>$. \et

\begin{proof} This had been essentially proved in \cite{li-qva1}
though the notion of nonlocal vertex $\F((t))$-algebra was absent.
It was proved in \cite{li-qva1} (Proposition 3.12) that if
$(a(x),b(x))$ is a quasi compatible ordered pair in $\E(W)$, then
for any $f(x),g(x)\in \F((x))$, $(f(x)a(x),g(x)b(x))$ is a quasi
compatible ordered pair and
\begin{eqnarray}\label{eysemi-relation}
Y_{\E}(f(x)a(x),x_{0})(g(x)b(x))=f(x+x_{0})g(x)Y_{\E}(a(x),x_{0})b(x).
\end{eqnarray}
It was also proved that the $\F((x))$-span of any $Y_{\E}$-closed
quasi compatible $\F$-subspace of $\E(W)$ is quasi compatible and
$Y_{\E}$-closed. It can be readily seen from the proof that this is
also true for compatible case. The rest follows from Theorem
\ref{tquasi-main}.
\end{proof}

We continue to establish a construction of weak quantum vertex
$\F((t))$-algebras.

\bd{dpseudo-local} {\em Let $W$ be a vector space over $\F$ as
before. A subset $U$ of $\E(W)$ is said to be {\em quasi
$\S(x_{1},x_{2})$-local}  if for any $a(x),b(x)\in U$, there exist
finitely many
$$u^{(i)}(x),\ v^{(i)}(x)\in U,\ \; f_{i}(x_{1},x_{2})\in \F_{*}(x_{1},x_{2})
\ \ (i=1,\dots,r)$$ such that
\begin{eqnarray}\label{epseudo-local}
p(x_{1},x_{2})a(x_{1})b(x_{2})
=\sum_{i=1}^{r}p(x_{1},x_{2})\iota_{x_{2},x_{1}}(f_{i}(x_{1},x_{2}))
u^{(i)}(x_{2})v^{(i)}(x_{1})
\end{eqnarray}
for some nonzero $p(x_{1},x_{2})\in \F[x_{1},x_{2}]$, depending on
$a(x)$ and $b(x)$. We say that $U$ is {\em $\S(x_{1},x_{2})$-local}
if the polynomial $p(x_{1},x_{2})$ is of the form
$(x_{1}-x_{2})^{k}$ with $k\in \N$. } \ed

We note that a quasi $\S(x_{1},x_{2})$-local subset is the same as a
pseudo-local subset as defined in \cite{li-qva1}. The following is
straightforward to prove:

\bl{lspan-qcomp} The $\F((x))$-span of any (resp. quasi)
$\S(x_{1},x_{2})$-local subset of $\E(W)$ is (resp. quasi)
$\S(x_{1},x_{2})$-local. \el

We also have:

\bp{ppseudo} Every (resp. quasi) $\S(x_{1},x_{2})$-local subset $U$
of $\E(W)$ is (resp. quasi) compatible. Furthermore, the
$\F((x))$-submodule $\F((x))\<U\>$ as in Theorem
\ref{tconcrete-general} is (resp. quasi) $\S(x_{1},x_{2})$-local.
\ep

\begin{proof} It was proved in \cite{li-qva1} (Lemma 3.2 and
Proposition 3.9) that if $U$ is quasi $\S(x_{1},x_{2})$-local, $U$
is quasi compatible and $\<U\>$ is quasi $\S(x_{1},x_{2})$-local.
Following the same proof with the obvious changes, we confirm the
corresponding assertions without the word ``quasi'' in the three
places. Then, by Lemma \ref{lspan-qcomp} $\F((x))\<U\>$ is (resp.
quasi) $\S(x_{1},x_{2})$-local.
\end{proof}

Furthermore, we have:

\bp{padjoint-relation} Let $W$ be a vector space over $\F$ as before
and let $V$ be a $Y_{\E}$-closed quasi compatible
$\F((x))$-submodule of $\E(W)$. Let
$$a(x),\ b(x), \ u^{(i)}(x),\ v^{(i)}(x)\in V, \
f_{i}(x_{1},x_{2})\in \F_{*}(x_{1},x_{2}) \; \; (i=1,\dots,r).$$
Assume that there exists a nonzero polynomial $p(x_{1},x_{2})\in
\F[x_{1},x_{2}]$ such that
$$p(x_{1},x_{2})a(x_{1})b(x_{2})
=\sum_{i=1}^{r}p(x_{1},x_{2})\iota_{x_{2},x_{1}}(f_{i}(x_{1},x_{2}))
u^{(i)}(x_{2})v^{(i)}(x_{1}).$$ Then
\begin{eqnarray}
& &(x_{1}-x_{2})^{k}Y_{\E}(a(x),x_{1})Y_{\E}(b(x),x_{2})\nonumber\\
&=&(x_{1}-x_{2})^{k}
\sum_{i=1}^{r}\iota_{x,x_{2},x_{1}}(f_{i}(x+x_{1},x+x_{2}))
Y_{\E}(u^{(i)}(x),x_{2})Y_{\E}(v^{(i)}(x),x_{1}),\ \ \ \
\end{eqnarray}
where $p(x_{1},x_{2})=(x_{1}-x_{2})^{k}q(x_{1},x_{2})$ with $k\in
\N,\; q(x_{1},x_{2})\in \F[x_{1},x_{2}]$ such that
$q(x_{1},x_{1})\ne 0$. \ep

\begin{proof} By Proposition 3.13 of \cite{li-qva1} we have
\begin{eqnarray*}
& &p(x+x_{1},x+x_{2})Y_{\E}(a(x),x_{1})Y_{\E}(b(x),x_{2})\nonumber\\
&=&p(x+x_{1},x+x_{2})
\sum_{i=1}^{r}\iota_{x,x_{2},x_{1}}(f_{i}(x+x_{1},x+x_{2}))
Y_{\E}(u^{(i)}(x),x_{2})Y_{\E}(v^{(i)}(x),x_{1}).
\end{eqnarray*}
Note that
$$p(x+x_{1},x+x_{2})=(x_{1}-x_{2})^{k}q(x+x_{1},x+x_{2}).$$
 Write $q(x+x_{1},x+x_{2})=q(x,x)+x_{1}g+x_{2}h$ with
$g,h\in \F[x,x_{1},x_{2}]$. As $q(x,x)\ne 0$, by Lemma
\ref{lsubcomm} we have
$$\iota_{x,x_{1},x_{2}}(q(x+x_{1},x+x_{2})^{-1})
=\iota_{x,x_{2},x_{1}}(q(x+x_{1},x+x_{2})^{-1})\in
\F((x))[[x_{1},x_{2}]].$$ Then we can cancel the factor
$q(x+x_{1},x+x_{2})$ to obtain the desired relation.
\end{proof}

As our second main result of this section we have the following
refinement of Theorem \ref{tconcrete-general}:

\bt{tmain} Let $W$ be a vector space over $\F$ and let $U$ be any
(resp. quasi) $\S(x_{1},x_{2})$-local subset of $\E(W)$. Then the
nonlocal vertex $\F((t))$-algebra $\F((x))\<U\>$ which was obtained
in Theorem \ref{tconcrete-general} is a weak quantum vertex
$\F((t))$-algebra with $W$ as a type zero (resp. quasi) module.\et

\begin{proof} Since a (resp. quasi) $\S(x_{1},x_{2})$-local subset
is (resp. quasi) compatible by Proposition \ref{ppseudo}, the
assertion on module structure follows from Theorem
\ref{tconcrete-general}. As for the first assertion, by Proposition
\ref{ppseudo}, $\F((x))\<U\>$ is (resp. quasi)
$\S(x_{1},x_{2})$-local. Then it follows from Proposition
\ref{padjoint-relation} that the nonlocal vertex $\F((t))$-algebra
$\F((x))\<U\>$ satisfies $\S_{t}$-locality. In view of Proposition
\ref{ptlocality}, $\F((x))\<U\>$ is a weak quantum vertex
$\F((t))$-algebra.
\end{proof}

We end this section with the following technical result:

\bp{prelation-am} Let $W$ be a vector space over $\F$, let $V$ be a
$Y_{\E}$-closed quasi compatible $\F((x))$-submodule of $\E(W)$, and
let
$$n\in \Z,\ a(x),\ b(x),\ u^{(i)}(x),\ v^{(i)}(x)\in V,\ f_{i}(x_{1},x_{2})\in
\F_{*}(x_{1},x_{2}) \ \ (i=1,\dots,r),$$
$$c^{(0)}(x),c^{(1)}(x),\dots,c^{(s)}(x)\in V.$$
Assume
\begin{eqnarray}\label{ecross-product}
& &(x_{1}-x_{2})^{n}a(x_{1})b(x_{2})-(-x_{2}+x_{1})^{n}
\sum_{i=1}^{r}\iota_{x_{2},x_{1}}(f_{i}(x_{1},x_{2}))
u^{(i)}(x_{2})v^{(i)}(x_{1})\nonumber\\
&=&\sum_{j=0}^{s}c^{(j)}(x_{2})\frac{1}{j!}\left(\frac{\partial}{\partial
x_{2}}\right)^{j} x_{1}^{-1}\delta\left(\frac{x_{2}}{x_{1}}\right).
\end{eqnarray}
Then
\begin{eqnarray}\label{ecross-product-va}
& &(x_{1}-x_{2})^{n}Y_{\E}(a(x),x_{1})Y_{\E}(b(x),x_{2})\nonumber\\
& &\ \ \ \ \ -(-x_{2}+x_{1})^{n}
\sum_{i=1}^{r}\iota_{x,x_{2},x_{1}}(f_{i}(x+x_{1},x+x_{2}))
Y_{\E}(u^{(i)}(x),x_{2})Y_{\E}(v^{(i)}(x),x_{1})\nonumber\\
&=&\sum_{j=0}^{s}Y_{\E}(c^{(j)}(x),x_{2})\frac{1}{j!}\left(\frac{\partial}{\partial
x_{2}}\right)^{j} x_{1}^{-1}\delta\left(\frac{x_{2}}{x_{1}}\right).
\end{eqnarray}
 \ep

\begin{proof} Let $k$ be a nonnegative integer such that $k\ge s+1$ and $k+n\ge 0$.
Multiplying both sides of (\ref{ecross-product}) by
$(x_{1}-x_{2})^{k}$ we obtain
\begin{eqnarray}\label{eproof-slocal}
(x_{1}-x_{2})^{k+n}a(x_{1})b(x_{2}) =(x_{1}-x_{2})^{k+n}
\sum_{i=1}^{r}\iota_{x_{2},x_{1}}(f_{i}(x_{1},x_{2}))
u^{(i)}(x_{2})v^{(i)}(x_{1}),
\end{eqnarray}
as $(x_{1}-x_{2})^{k}\left(\frac{\partial}{\partial
x_{2}}\right)^{j}
x_{1}^{-1}\delta\left(\frac{x_{2}}{x_{1}}\right)=0$ for $0\le j\le
s$. By Proposition \ref{padjoint-relation}, we have
\begin{eqnarray*}
& &(x_{1}-x_{2})^{n+k}Y_{\E}(a(x),x_{1})Y_{\E}(b(x),x_{2})\nonumber\\
&=&(x_{1}-x_{2})^{n+k}\sum_{i=1}^{r}\iota_{t,x_{2},x_{1}}(f_{i}(t+x_{1},t+x_{2}))
Y_{\E}(u^{(i)}(x),x_{2})Y_{\E}(v^{(i)}(x),x_{1}),\ \  \  \
\end{eqnarray*}
which together with weak associativity implies (by Lemma
\ref{lformal-jacobi-old})
\begin{eqnarray}\label{eproof-jacobi}
& &x_{0}^{-1}\delta\left(\frac{x_{1}-x}{x_{0}}\right)
Y_{\E}(a(x),x_{1})Y_{\E}(b(x),x_{2})\nonumber\\
& & \ \
-x_{0}^{-1}\delta\left(\frac{x-x_{1}}{-x_{0}}\right)\sum_{i=1}^{r}
\iota_{t,x_{2},x_{1}}(f_{i}(t+x_{1},t+x_{2}))
Y_{\E}(u^{(i)}(x),x_{2})Y_{\E}(v^{(i)}(x),x_{1})\nonumber\\
&=&x_{1}^{-1}\delta\left(\frac{x_{2}+x_{0}}{x_{1}}\right)
Y_{\E}(Y_{\E}(a(x),x_{0})b(x),x_{2}).
\end{eqnarray}
With (\ref{eproof-slocal}) we have
\begin{eqnarray*}
& &Y_{\E}(a(x),x_{0})b(x)=x_{0}^{-n-k}\Res_{x_{1}}
x_{1}^{-1}\delta\left(\frac{x+x_{0}}{x_{1}}\right)\left((x_{1}-x)^{n+k}a(x_{1})b(x)\right) \\
&&\ \ \ =\Res_{x_{1}}x_{0}^{-1}\delta\left(\frac{x_{1}-x}{x_{0}}\right)a(x_{1})b(x)\\
& &\hspace{1cm}
-\Res_{x_{1}}x_{0}^{-1}\delta\left(\frac{x-x_{1}}{-x_{0}}\right)\sum_{i=1}^{r}
\iota_{x,x_{1}}(f_{i}(x_{1},x))u^{(i)}(x)v^{(i)}(x_{1}),
\end{eqnarray*}
{}from which we obtain
\begin{eqnarray*}
a(x)_{n}b(x)=c^{(0)}(x),\ \ a(x)_{n+1}b(x)=c^{(1)}(x),\dots,
a(x)_{n+s}b(x)=c^{(s)}(x),
\end{eqnarray*}
and $a(x)_{m}b(x)=0$ for $m>n+s$. Then applying
$\Res_{x_{0}}x_{0}^{n}$ to (\ref{eproof-jacobi}) we obtain
(\ref{ecross-product-va}).
\end{proof}

\section{General existence theorems}
In this section we present two existence theorems for a nonlocal
vertex $\F((t))$-algebra structure and for a weak quantum vertex
$\F((t))$-algebra structure. These are analogs of the existence
theorem in the theory of weak quantum vertex algebras (see
\cite{li-qva1}, \cite{li-qva2}) and in the theory of vertex algebras
(see \cite{fkrw}, \cite{mp}; cf. \cite{ll}).

We begin by reexamining Section 4 with $\F((t))$ in place of $\F$ as
the scalar field. Let $W$ be an $\F((t))$-module, namely a vector
space over $\F((t))$. By $\E(W)$ we mean the $\F((t))$-module
\begin{eqnarray}
\E(W)=\Hom_{\F((t))}(W,W((x))),
\end{eqnarray}
which is a canonical $\F((t))((x))$-module. Let $W_{\F}$ denote $W$
viewed as a vector space over $\F$. We see that $\E(W)\subset
\E(W_{\F})$ and that every compatible subset of $\E(W)$ is also a
compatible subset of $\E(W_{\F})$.

As a convention, for $f(t)\in \F((t))$  we define
$$f(t+x)=\iota_{t,x}f(t+x)=e^{x\frac{d}{dt}}f(t)\in \F((t))[[x]]\subset \F((t))((x)).$$

The following is immediate:

\bl{lt1module} Let $W$ be an $\F((t))$-module and let $t_{1}$ be
another formal variable. Then $\E(W)$ becomes an
$\F((t_{1}))$-module with
\begin{eqnarray} f(t_{1})a(x)=f(t+x)a(x)\
\ \ \mbox{ for }f(t_{1})\in \F((t_{1})),\; a(x)\in \E(W).
\end{eqnarray}
\el

With this we have:

\bp{pcompatible-t} Let $W$ be an $\F((t))$-module and let $U$ be a
compatible subset of $\E(W)$. Denote by $\<U\>_{\F}$ the nonlocal
vertex algebra over $\F$ generated by $U$. Then
$\F((t_{1}))\<U\>_{\F}$ is a nonlocal vertex $\F((t_{1}))$-algebra,
and $W$, viewed as an $\F((t_{1}))$-module with $f(t_{1})\in
\F((t_{1}))$ acting as $f(t)$, is a module of type one. \ep

\begin{proof} By Theorem \ref{tquasi-main} with $\F((t))$ in place of $\F$,
$U$ generates a nonlocal vertex algebra $\<U\>$ over $\F((t))$.
Furthermore, by Theorem \ref{tconcrete-general} the span
$\F((t))((x))\<U\>$ is also a nonlocal vertex algebra over
$\F((t))$, satisfying the condition that
$$Y_{\E}(g(x)a(x),x_{0})(h(x)b(x))=g(x+x_{0})h(x)Y_{\E}(a(x),x_{0})b(x)$$
for $g(x),h(x)\in \F((t))((x)),\; a(x),b(x)\in \F((t))((x))\<U\>$.
{}From this we have
\begin{eqnarray*}
Y_{\E}(f(t_{1})a(x),x_{0})(g(t_{1})b(x))&=&Y_{\E}(f(t+x)a(x),x_{0})(g(t+x)b(x))\\
&=&f(t+x+x_{0})g(t+x)Y_{\E}(a(x),x_{0})b(x)\\
&=&f(t_{1}+x_{0})g(t_{1})Y_{\E}(a(x),x_{0})b(x)
\end{eqnarray*}
for $f(t_{1}),g(t_{1})\in \F((t_{1})),\; a(x),b(x)\in
\F((t))((x))\<U\>$.  It follows that $\F((t))((x))\<U\>$ is a
nonlocal vertex $\F((t_{1}))$-algebra with $\F((t_{1}))\<U\>_{\F}$
as a subalgebra. Also, by Theorem \ref{tquasi-main}, $W$ is a module
for $\F((t))((x))\<U\>$ viewed as a nonlocal vertex algebra over
$\F((t))$ with $Y_{W}(\alpha(x),x_{0})=\alpha(x_{0})$ for
$\alpha(x)\in \F((t))((x))\<U\>$. For $f(t_{1})\in \F((t_{1})),\;
a(x)\in \F((t))((x))\<U\>,\; w\in W$, we have
\begin{eqnarray*}
Y_{W}(f(t_{1})a(x),x_{0})w&=&Y_{W}(f(t+x)a(x),x_{0})w=f(t+x_{0})a(x_{0})w\\
&=&f(t_{1}+x_{0})Y_{W}(a(x),x_{0})w.
\end{eqnarray*}
Then the last assertion follows.
\end{proof}

\bd{dstlocal} {\em Let $W$ be an $\F((t))$-module. A subset $U$ of
$\E(W)$ is said to be {\em $\S_{t}$-local} if for any $a(x),b(x)\in
U$, there exist (finitely many)
$$u^{(i)}(x),v^{(i)}(x)\in U,\; f_{i}(x_{1},x_{2})\in
\F_{*}(x_{1},x_{2})\;\; (i=1,\dots, r)$$ such that
\begin{eqnarray}\label{estlocal}
(x_{1}-x_{2})^{k}a(x_{1})b(x_{2})
=(x_{1}-x_{2})^{k}\sum_{i=1}^{r}\iota_{t,x_{2},x_{1}}(f_{i}(t+x_{1},t+x_{2}))
u^{(i)}(x_{2})v^{(i)}(x_{1})
\end{eqnarray}
for some nonnegative integer $k$ depending on $a(x)$ and $b(x)$. }
\ed

With this notion we have:

\bt{tst1} Let $W$ be an $\F((t))$-module and let $U$ be an
$\S_{t}$-local subset of $\E(W)$. Then $U$ is compatible.
Furthermore, $U$ is an $\S_{t_{1}}$-local subset of the nonlocal
vertex $\F((t_{1}))$-algebra $\F((t_{1}))\<U\>_{\F}$, which was
obtained in Proposition \ref{pcompatible-t}, and
$\F((t_{1}))\<U\>_{\F}$ is a weak quantum vertex
$\F((t_{1}))$-algebra and $W$, viewed as an $\F((t_{1}))$-module
with $f(t_{1})\in \F((t_{1}))$ acting as $f(t)$, is a module of type
one. \et

\begin{proof} Let $f(x_{1},x_{2})\in \F_{*}(x_{1},x_{2})$. We have
$f(t+x_{1},t+x_{2})\in \F_{*}(t,x_{1},x_{2})$ with
$$\iota_{t,x_{2},x_{1}}f(t+x_{1},t+x_{2})\in \F((t))((x_{2}))((x_{1})).$$
We can also view $f(t+x_{1},t+x_{2})$ as an element of
$\F((t))_{*}(x_{1},x_{2})$ (with $\F((t))$ as the scalar field),
which we denote by $f_{t}(x_{1},x_{2})$, noticing that for $q/p\in
\F_{*}(x_{1},x_{2})$ with $q\in \F[[x_{1},x_{2}]],\; p\in
\F[x_{1},x_{2}]$, we have
\begin{eqnarray*}
&&q(t+x_{1}, t+x_{2})\in \F[[t,x_{1},x_{2}]]\subset
\F((t))[[x_{1},x_{2}]],\\
&&p(t+x_{1},t+x_{2})\in \F[t,x_{1},x_{2}]\subset
\F((t))[x_{1},x_{2}].
\end{eqnarray*}
With the iota-map $\iota_{x_{2},x_{1}}:
\F((t))_{*}(x_{1},x_{2})\rightarrow \F((t))((x_{2}))((x_{1}))$, we
have
\begin{eqnarray*}
\iota_{x_{2},x_{1}}f_{t}(x_{1},x_{2})\in \F((t))((x_{2}))((x_{1})).
\end{eqnarray*}
It is straightforward to show that
$$\iota_{x_{2},x_{1}}f_{t}(x_{1},x_{2})=\iota_{t,x_{2},x_{1}}f(t+x_{1},t+x_{2}).$$
In view of this, we see that an $\S_{t}$-local subset of $\E(W)$ is
also an $\S(x_{1},x_{2})$-local subset with $\F((t))$ in place of
$\F$. By Proposition \ref{ppseudo}, every $\S_{t}$-local subset is
compatible.

Let $a(x),b(x)\in U$. There exist
$$u^{(i)}(x),v^{(i)}(x)\in U,\; f_{i}(x_{1},x_{2})\in
\F_{*}(x_{1},x_{2})\;\; (i=1,\dots, r)$$ such that (\ref{estlocal})
holds for some nonnegative integer $k$. Viewing
$f_{i}(t+x_{1},t+x_{2})$ as elements of $\F((t))_{*}(x_{1},x_{2})$,
from Proposition \ref{prelation-am} with $\F((t))$ in place of $\F$,
we have
\begin{eqnarray*}
&&(x_{1}-x_{2})^{k}Y_{\E}(a(x),x_{1})Y_{\E}(b(x),x_{2})\\
&=&(x_{1}-x_{2})^{k}\sum_{i=1}^{r}\iota_{t,x,x_{2},x_{1}}(f_{i}(t+x+x_{1},t+x+x_{2}))
Y_{\E}(u^{(i)}(x),x_{2})Y_{\E}(v^{(i)}(x),x_{1})\\
&=&(x_{1}-x_{2})^{k}\sum_{i=1}^{r}\iota_{t_{1},x_{2},x_{1}}(f_{i}(t_{1}+x_{1},t_{1}+x_{2}))
Y_{\E}(u^{(i)}(x),x_{2})Y_{\E}(v^{(i)}(x),x_{1}).
\end{eqnarray*}
This proves that $U$ is an $\S_{t_{1}}$-local subset of the nonlocal
vertex $\F((t_{1}))$-algebra $\F((t_{1}))\<U\>_{\F}$. Then by Lemma
\ref{lgenerating-qva}, $\F((t_{1}))\<U\>_{\F}$ is a weak quantum
vertex $\F((t_{1}))$-algebra. The last assertion on module structure
has already been established in Proposition \ref{pcompatible-t}.
\end{proof}

Now, we are ready to present our first existence theorem.

\bt{tconstruction} Let $V$ be an $\F((t))$-module, ${\bf 1}$ a
vector of $V$, $\D$ an $\F$-linear operator on $V$, $U$ an
$\F((t))$-submodule of $V$,
\begin{eqnarray*}
Y_{0}(\cdot,x):&& U\rightarrow \E(V)=\Hom_{\F((t))}(V,V((x)))\\
&&u\mapsto Y_{0}(u,x)=u(x)=\sum_{n\in \Z}u_{n}x^{-n-1},
\end{eqnarray*}
an $\F$-linear map, satisfying all the following conditions: $\D{\bf
1}=0$,
\begin{eqnarray*}
&&\D(f(t)v)=f(t)\D v+f'(t)v\ \ \ \mbox{ for }f(t)\in \F((t)),\; v\in V,\\
& &\hspace{1cm}[\D,Y_{0}(u,x)]=\frac{d}{dx}Y_{0}(u,x),\\
&&\hspace{1cm}Y_{0}(u,x){\bf 1}\in V[[x]]\ \mbox{ and }\
\lim_{x\rightarrow 0}Y_{0}(u,x){\bf 1}=u,\\
&&Y_{0}( f(t)u,x)=f(t+x)Y_{0}(u,x) \ \ \mbox{ for }u\in U,\; f(t)\in
\F((t)),
\end{eqnarray*}
$\{Y_{0}(u,x)\;|\; u\in U\}$, denoted by $U(x)$, is compatible, and
$V$ is linearly spanned over $\F((t))$ by the vectors
$$u^{(1)}_{m_{1}}\cdots u^{(r)}_{m_{r}}{\bf 1}$$
for $r\ge 0,\; u^{(i)}\in U,\; m_{i}\in \Z$. Suppose that there
exists an $\F$-linear map $\psi_{x}$ from $V$ to
$\F((t))((x))\<U(x)\>$ such that
\begin{eqnarray*}
&&\psi_{x}(f(t)v)=f(t+x)\psi_{x}(v)\ \ \
\mbox{ for }f(t)\in \F((t)),\; v\in V,\\
&&\psi_{x}({\bf 1})=1_{V},\ \ \
\psi_{x}(u_{n}v)=u(x)_{n}\psi_{x}(v)\ \ \mbox{ for }u\in U,\; n\in
\Z,\; v\in V.
\end{eqnarray*}
For $v\in V$, set $Y(v,x)=\psi_{x}(v)\in \E(V)$. Then $Y(\cdot,x)$
extends $Y_{0}(\cdot,x)$, and $(V,Y,{\bf 1})$ carries the structure
of a nonlocal vertex $\F((t))$-algebra. Furthermore, if $U(x)$ is
$\S_{t}$-local, $V$ is a weak quantum vertex $\F((t))$-algebra. \et

\begin{proof} First consider the case with ${\bf 1}\in U$.
Then $\F((t)){\bf 1}\subset U$.
 Since $(f(t){\bf 1})_{-1}{\bf
1}=f(t){\bf 1}$ for $f(t)\in \F((t))$, we see that $V$ is actually
linearly spanned over $\F$ by those vectors in the assumption. It
follows from \cite{li-qva1} (Theorem 6.3) that $(V,Y,{\bf 1})$
carries the structure of a nonlocal vertex algebra over $\F$. It is
clear that $U$ generates $V$ as a nonlocal vertex algebra over $\F$
and we have
$$Y(f(t)u,x)g(t)=Y_{0}(f(t)u,x)g(t)=f(t+x)g(t)Y_{0}(u,x)=f(t+x)g(t)Y(u,x)$$
for $f(t),g(t)\in \F((t)),\; u\in U$.  It follows from Lemma
\ref{lgenerating-1} that $V$ is a nonlocal vertex $\F((t))$-algebra.
Furthermore, if $\{Y_{0}(u,x)\;|\; u\in U\}$ ($=U(x)$) is
$\S_{t}$-local, it follows from Lemma \ref{lgenerating-qva} that $V$
is a weak quantum vertex $\F((t))$-algebra.

Now, assume ${\bf 1}\notin U$. Then $U\cap \F((t)){\bf 1}=0$. Set
$\bar{U}=U\oplus\F((t)){\bf 1}$. Extend the map $Y_{0}$ to $\bar{U}$
by defining $\bar{Y}_{0}(f(t){\bf 1},x)=f(t+x)$ for $f(t)\in
\F((t))$. We have
$$[\D,\bar{Y}_{0}(f(t){\bf 1},x)]=[\D,f(t+x)]=f'(t+x)
=\frac{d}{dx}\bar{Y}_{0}(f(t){\bf 1},x),$$
$$\bar{Y}_{0}(f(t){\bf 1},x){\bf 1}=f(t+x){\bf 1}\in V[[x]]\ \
\mbox{ and }\ \ \lim_{x\rightarrow 0}\bar{Y}_{0}(f(t){\bf 1},x){\bf
1}=f(t){\bf 1}.$$ Noticing that for $f(t)\in \F((t)),\; a(x)\in
\E(V)$,
$$Y_{\E}(f(t+x),x_{0})a(x)=(f(t+x_{1})a(x))|_{x_{1}=x+x_{0}}=f(t+x+x_{0})a(x),$$
we get $$\psi_{x}(\bar{Y}_{0}(f(t){\bf
1},x_{0})v)=\psi_{x}(f(t+x_{0})v)=f(t+x+x_{0})\psi_{x}(v)=Y_{\E}(f(t+x),x_{0})\psi_{x}(v)$$
for $v\in V$. Then it follows from the first part with $(U,Y_{0})$
in place of $(\bar{U},\bar{Y}_{0})$.
\end{proof}

In fact, for the last assertion of Theorem \ref{tconstruction} on
the weak quantum vertex $\F((t))$-algebra structure, we can remove
those assumptions involving the operator $\D$ (cf. \cite{li-qva2},
Theorem 2.9). First we prove the following (cf. \cite{li-qva2},
Proposition 2.8):

\bl{lgenerating-3} Let $V$ be a nonlocal vertex $\F((t))$-algebra
and let $(W,Y_{W})$ be a type one $V$-module.  Suppose that $e$ is a
vector in $W$ and $U$ is an $\S_{t}$-local subset of $V$ such that
$$Y_{W}(u,x)e\in W[[x]]\ \ \ \mbox{ for }u\in U.$$
Then $Y_{W}(v,x)e\in W[[x]]$ for $v\in \F((t))\<U\>_{\F}$ and the
map $\theta: \F((t))\<U\>_{\F}\rightarrow W,$ sending $v$ to
$v_{-1}e$ for $v\in \F((t))\<U\>_{\F}$, is a homomorphism of type
one $\F((t))\<U\>_{\F}$-modules  with $\theta({\bf 1})=e$. \el

\begin{proof}  Set
$$K=\{ v\in V\; |\; Y_{W}(v,x)e\in W[[x]]\}.$$
It can be readily seen that $K$ is an $\F((t))$-submodule of $V$. We
must prove $\F((t))\<U\>_{\F}\subset K$. {}From assumption,
$\F((t)){\bf 1}+\F((t))U$ is an $\S_{t}$-local $\F((t))$-submodule
of $K$. Then there exists a maximal $\S_{t}$-local
$\F((t))$-submodule $A$ of $K$, containing
 $\F((t)){\bf 1}+\F((t))U$. We now prove $\F((t))\<U\>_{\F}\subset A\ (\subset K)$.
 As $\{{\bf 1}\}\cup
U\subset A$, it suffices to prove that $A$ is closed, i.e.,
$A^{(2)}\subset A$.  {}From the proof of Lemma
\ref{lgenerating-qva}, $A^{(2)}$ is $\S_{t}$-local. Now we prove
$A^{(2)}\subset K$. Let $u,v\in A$. From Proposition
\ref{pmodule-property-type2}, there exist
$$u^{(i)},\ v^{(i)}\in A,\; f_{i}(x_{1},x_{2})\in
\F_{*}(x_{1},x_{2}) \ \ (i=1,\dots,r)$$ such that the Jacobi
identity (\ref{et-jacobi-module-type2}) holds. By Lemma
\ref{lsubcomm}, for $1\le i\le r$, the series
$$\iota_{t,x_{2},x_{1}}(f_{i}(t+x_{1},t+x_{2}))$$ involves only
nonnegative powers of $x_{1}$.  By applying $\Res_{x_{1}}$ to
(\ref{et-jacobi-module-type2}), we see that
$$Y_{W}(Y(u,x_{0})v,x_{2})e\in W[[x_{2}]]((x_{0})),$$ which implies
$u_{m}v\in K$ for $m\in \Z$. Thus $A^{(2)}\subset K$. Since ${\bf
1}\in A$, we have $A\subset A^{(2)}$. As $A$ is maximal we must have
$A=A^{(2)}$, proving that $A$ is closed. Thus we have
$\F((t))\<U\>_{\F}\subset A\subset K$, proving the first assertion.

By Lemma 6.1 of \cite{li-qva1}, $\theta$ is a module homomorphism
for $\F((t))\<U\>_{\F}$ viewed as a nonlocal vertex algebra over
$\F$. Furthermore, from Lemma \ref{lconsistency}, $\theta$ is
$\F((t))$-linear. Thus $\theta$ is a homomorphism of type one
$\F((t))\<U\>_{\F}$-modules.
\end{proof}

Now, we have:

\bt{tconstruction-refine} Let $V, {\bf 1}, U, Y_{0}(\cdot,x), U(x)$,
and $\psi_{x}$ be given as in Theorem \ref{tconstruction} and retain
all the assumptions that do not involve $\D$. In addition, assume
that $U(x)$ is $\S_{t}$-local. Set $Y(v,x)=\psi_{x}(v)\in \E(V)$ for
$v\in V$. Then $Y(\cdot,x)$ extends the map $Y_{0}(\cdot,x)$, and
$(V,Y,{\bf 1})$ carries the structure of a weak quantum vertex
$\F((t))$-algebra. \et

\begin{proof} Recall from Lemma \ref{lt1module} the
$\F((t_{1}))$-module structure on  $\E(V)$ with $t_{1}$ another
formal variable. For $u\in U,\; f(t_{1})\in \F((t_{1}))$, we have
$$f(t_{1})u(x)=f(t+x)Y_{0}(u,x)=Y_{0}(f(t)u,x)=(f(t)u)(x),$$
so that $U(x)$ is an $\F((t_{1}))$-submodule of $\E(V)$. Since
$U(x)$ is an $\S_{t}$-local subset of $\E(V)$, by Theorem
\ref{tst1}, $U(x)$ is compatible and $\F((t_{1}))\<U(x)\>_{\F}$ is a
weak quantum $\F((t_{1}))$-algebra, where $\<U(x)\>_{\F}$ denotes
the nonlocal vertex algebra over $\F$, generated by $U(x)$ inside
$\E(V)$. Set $E=\F((t_{1}))\<U(x)\>_{\F}$. By Theorem \ref{tst1},
$U(x)$ is an $\S_{t_{1}}$-local $\F((t_{1}))$-submodule of $E$ and
$(V,Y_{V})$ is a type one $E$-module, where $V$ is viewed as an
$\F((t_{1}))$-module with $f(t_{1})$ acting as $f(t)$ and where
$Y_{V}(a(x),x_{0})=a(x_{0})$ for $a(x)\in E.$ {}From our assumption
we have
$$Y_{V}(u(x),x_{0}){\bf 1}=u(x_{0}){\bf 1}\in V[[x_{0}]]
\ \ \mbox{ for }u\in U. $$ By Lemma \ref{lgenerating-3}, there
exists an $E$-module homomorphism $\phi$ from $E$ to $V$ such that
$\phi(1_{V})={\bf 1}$ and
$$\phi(f(t+x)a(x))=\phi(f(t_{1})a(x))=f(t_{1})\phi(a(x))=f(t)\phi(a(x))$$
for $f(t)\in \F((t)),\; a(x)\in E$. For $u\in U,\; a(x)\in E$, we
have
$$\phi(Y_{\E}(u(x),x_{0})a(x))=Y_{V}(u(x),x_{0})\phi(a(x))=u(x_{0})\phi(a(x)).$$
That is,
$$\phi(u(x)_{n}a(x))=u_{n}\phi(a(x))\ \ \mbox{ for } u\in
U,\; a(x)\in E,\; n\in \Z.$$ It follows that $\psi_{x}$ is an
$\F((t_{1}))$-isomorphism from $V$ to $E$ with $\phi$ as the
inverse. Then we have a weak quantum vertex $\F((t_{1}))$-algebra
structure on $V$, transported from $E$, where ${\bf 1}\
(=\phi(1_{V}))$ is the vacuum vector. The defined map $Y(\cdot,x)$
coincides with the transported structure, as for $v\in V$,
$$\phi Y_{\E}(\psi_{x}(v),x_{0})\psi_{x}
=Y_{V}(\psi_{x}(v),x_{0})=\psi_{x_{0}}(v)=Y(v,x_{0})$$ (recall that
$\phi$ is a module homomorphism). Furthermore, for $u\in U$, we have
$$\psi_{x}(u)=\psi_{x}(u_{-1}{\bf 1})=u(x)_{-1}1_{V}=u(x),$$
so that
$$Y(u,x)=\psi_{x}(u)=u(x)=Y_{0}(u,x).$$
Now, the proof is completed.
\end{proof}

\section{Example of quantum vertex $\C((t))$-algebras}
In this section we first associate weak quantum $\C((t))$-algebras
to quantum affine algebras and then we construct an example of
non-degenerate quantum vertex $\C((t))$-algebras from a certain
quantum $\beta\gamma$-system.

First, we follow \cite{fj} (see also \cite{drqaa}) to present the
quantum affine algebras. Let $\g$ be a finite-dimensional simple Lie
algebra of rank $l$ of type $A$, $D$, or $E$ and let $A=(a_{ij})$ be
the Cartan matrix.  Let $q$ be a nonzero complex number. For $1\le
i,j\le l$, set
\begin{eqnarray}
f_{ij}(x)=(q^{a_{ij}}x-1)/(x-q^{a_{ij}})\in \C(x).
\end{eqnarray}
Furthermore, set
\begin{eqnarray}
g_{ij}(x)^{\pm 1}=\iota_{x,0} f_{ij}(x)^{\pm 1}\in \C[[x]],
\end{eqnarray}
where $\iota_{x,0} f_{ij}(x)^{\pm 1}$ are the formal Taylor series
expansions of $f_{ij}(x)^{\pm 1}$ at $0$.  The quantum affine
algebra $U_{q}(\hat{\g})$ is (isomorphic to) the associative algebra
with identity $1$ and with generators
\begin{eqnarray}
X_{ik}^{\pm},\ \ \phi_{im},\ \  \psi_{in},\ \ \gamma^{1/2},\ \
\gamma^{-1/2}
\end{eqnarray}
for $1\le i\le l, \;k\in \Z,\; m\in -\N,\;n\in \N$, where
$\gamma^{\pm 1/2}$ are central elements, satisfying the relations
below, written in terms of the following generating functions:
\begin{eqnarray}
X_{i}^{\pm}(z)=\sum_{k\in \Z}X_{ik}^{\pm}z^{-k},\ \ \ \
\phi_{i}(z)=\sum_{m\in -\N}\phi_{im}z^{-m},\ \ \ \
\psi_{i}(z)=\sum_{n\in \N}\psi_{in}z^{-n}.
\end{eqnarray}
The relations are
\begin{eqnarray*}
& &\gamma^{1/2}\gamma^{-1/2}=\gamma^{-1/2}\gamma^{1/2}=1,\\
& &\phi_{i0}\psi_{i0}=\psi_{i0}\phi_{i0}=1,\\
& &[\phi_{i}(z),\phi_{j}(w)]=0,\ \ \ \
[\psi_{i}(z),\psi_{j}(w)]=0,\\
& &\phi_{i}(z)\psi_{j}(w)\phi_{i}(z)^{-1}\psi_{j}(w)^{-1}
=g_{ij}(z/w\gamma)/g_{ij}(z\gamma/w),\\
& &\phi_{i}(z)X^{\pm}_{j}(w)\phi_{i}(z)^{-1}
=g_{ij}(z/w\gamma^{\pm 1/2})^{\pm 1}X^{\pm}_{j}(w),\\
& &\psi_{i}(z)X^{\pm}_{j}(w)\psi_{i}(z)^{-1}
=g_{ij}(w/z\gamma^{\pm 1/2})^{\mp 1}X^{\pm}_{j}(w),\\
& &(z-q^{\pm a_{ij}}w)X^{\pm}_{i}(z)X^{\pm}_{j}(w)
=(q^{\pm a_{ij}}z-w)X^{\pm}_{j}(w)X^{\pm}_{i}(z),\\
& &[X^{+}_{i}(z),X^{-}_{j}(w)]=\frac{\delta_{ij}}{q-q^{-1}}
\left(\delta\left(\frac{z}{w\gamma}\right)\psi_{i}(w\gamma^{1/2})
-\delta\left(\frac{z\gamma}{w}\right)\phi_{i}(z\gamma^{1/2})\right),
\ \ \ \
\end{eqnarray*}
and there is one more set of relations of Serre type.

A $U_{q}(\hat{\g})$-module $W$ is said to be {\em restricted} if
 for any $w\in W$, $X_{ik}^{\pm}w=0$ and
$\psi_{ik}w=0$ for $1\le i\le l$ and for $k$ sufficiently large. We
say $W$ is of {\em level} $\ell\in \C$ if $\gamma^{\pm 1/2}$ act on
$W$ as scalars $q^{\pm \ell/4}$. (Rigorously speaking, one needs to
choose a branch of $\log q$.) We have (cf. \cite{li-qva1},
Proposition 4.9):

\bp{pqaffine} Let $q$ and $\ell$ be complex numbers with $q\ne 0$
and let $W$ be a restricted $U_{q}(\hat{\g})$-module of level
$\ell$. Set
$$U_{W}=\{ \phi_{i}(x), \psi_{i}(x), X_{i}^{\pm}(x)
\;|\; 1\le i\le l\}.$$ Then $U_{W}$ is a quasi
$\S(x_{1},x_{2})$-local subset of $\E(W)$ and $\C((x))\<U_{W}\>$ is
a weak quantum vertex $\C((t))$-algebra with $W$ as a type zero
quasi module, where $\<U_{W}\>$ denotes the nonlocal vertex algebra
over $\C$ generated by $U_{W}$. \ep

\begin{proof} As $W$ is a restricted module, we have $U_{W}\subset \E(W)$, noticing that
$\phi_{i}(x)\in (\End W)[[x]]\subset \E(W)$.
 Note that
$$g_{ij}(z/w)=\iota_{w,z}(q^{a_{ij}}z-w)/(z-q^{a_{ij}}w).$$
Then
$$g_{ij}(z/w\gamma)/g_{ij}(z\gamma/w)
=\iota_{w,z}\frac{(q^{a_{ij}}z-w\gamma)(z\gamma-q^{a_{ij}}w)}{(z-q^{a_{ij}}w\gamma)(q^{a_{ij}}z\gamma-w)}.$$
With this, from the defining relation we get
\begin{eqnarray}
(z-q^{a_{ij}}w\gamma)(q^{a_{ij}}z\gamma-w)\phi_{i}(z)\psi_{j}(w)
=(q^{a_{ij}}z-w\gamma)(z\gamma-q^{a_{ij}}w)\psi_{j}(w)\phi_{i}(z).
\end{eqnarray}
With
$$g_{ij}(z/w\gamma^{\pm 1/2})
=\iota_{w,z}(q^{a_{ij}}z-w\gamma^{\pm
1/2})/(z-q^{a_{ij}}w\gamma^{\pm 1/2}),$$ we get
\begin{eqnarray}
&&(z-q^{a_{ij}}w\gamma^{1/2})\phi_{i}(z)X^{+}_{j}(w)
=(q^{a_{ij}}z-w\gamma^{1/2})X^{+}_{j}(w)\phi_{i}(z),\\
&&(q^{a_{ij}}z-w\gamma^{-1/2})\phi_{i}(z)X^{-}_{j}(w)
=(z-q^{a_{ij}}w\gamma^{-1/2})X^{-}_{j}(w)\phi_{i}(z).
\end{eqnarray}
Similarly, we have
\begin{eqnarray}
& &(w-q^{a_{ij}}z\gamma^{1/2})\psi_{i}(z)X^{+}_{j}(w)
=(q^{a_{ij}}w-z\gamma^{1/2})X^{+}_{j}(w)\psi_{i}(z),\\
& &(q^{a_{ij}}z-w\gamma^{1/2})\psi_{i}(z)X^{-}_{j}(w)
=(z-q^{a_{ij}}w\gamma^{1/2})X^{-}_{j}(w)\psi_{i}(z).
\end{eqnarray}
As $(z-x)\delta(\frac{z}{x})=0$, we have
\begin{eqnarray}
(z-w\gamma)(z\gamma-w)X^{+}_{i}(z)X^{-}_{j}(w)
=(z-w\gamma)(z\gamma-w)X^{-}_{j}(w)X^{+}_{i}(z).
\end{eqnarray}
Now, it is clear that $U_{W}$ is a quasi $\S(x_{1},x_{2})$-local
subset of $\E(W)$. The rest follows immediately from Theorem
\ref{tmain}.
\end{proof}

For the rest of section, we construct a quantum vertex
$\C((t))$-algebra from a quantum $\beta\gamma$-system.

\bd{dqbetagamma} {\em Let $q$ be a nonzero complex number. Define
$A_{q}(\beta\gamma)$ to be the associative algebra over $\C$ with
generators $\beta_{n},\gamma_{n}\; (n\in \Z)$, subject to relations
\begin{eqnarray*}
&&\beta(x_{1})\beta(x_{2})
=\left(\frac{qx_{2}-x_{1}}{x_{2}-qx_{1}}\right)\beta(x_{2})\beta(x_{1}),\\
&&\gamma(x_{1})\gamma(x_{2})=\left(\frac{qx_{2}-x_{1}}{x_{2}-qx_{1}}\right)\gamma(x_{2})\gamma(x_{1}),\\
&&\beta(x_{1})\gamma(x_{2})-\left(\frac{x_{2}-qx_{1}}{qx_{2}-x_{1}}\right)\gamma(x_{2})\beta(x_{1})
=x_{1}^{-1}\delta\left(\frac{x_{2}}{x_{1}}\right),
\end{eqnarray*}
where $\beta(x)=\sum_{n\in \Z}\beta_{n}x^{-n-1}$ and
$\gamma(x)=\sum_{n\in \Z}\gamma_{n}x^{-n-1}$.} \ed

This algebra $A_{q}(\beta\gamma)$ belongs to a family of algebras,
known as Zamolodchikov-Faddeev algebras (see \cite{za}, \cite{fa}).
Notice that $A_{q}(\beta\gamma)$ becomes the standard
$\beta\gamma$-algebra when $q=1$, while $A_{q}(\beta\gamma)$ becomes
a Clifford algebra when $q=-1$. For these two special cases, it is
well known that a vertex algebra for $q=1$, or a vertex superalgebra
for $q=-1$ can be associated to the algebra $A_{q}(\beta\gamma)$
canonically. In the following, we shall mainly deal with the case
with $q \ne 1$.  (All the results will still hold for $q=1$, though
a different proof is needed.)

\br{rpreciseAq} {\em Notice that the defining relations involve
infinite sums, so that $A_{q}(\beta\gamma)$ is in fact a topological
algebra. One can give a precise definition using the procedure in
\cite{fz} for defining the universal enveloping algebra $U(V)$ of a
vertex operator algebra $V$. However, for this paper we shall only
need a category of modules for a free associative algebra. By an
{\em $A_{q}(\beta\gamma)$-module} we mean a vector space $W$ on
which the set $\{ \beta_{n},\gamma_{n}\;|\; n\in \Z\}$ acts as
linear operators, satisfying the condition that for every $w\in W$,
$$\beta_{n}w=\gamma_{n}w=0 \ \ \ \mbox{ for $n$ sufficiently large}$$
and all the relations in Definition \ref{dqbetagamma} after applied
to $w$ hold. } \er

\bd{dtbetagamma}{\em Let $q$ be a nonzero complex number as before.
Define $A_{t,q}(\beta\gamma)$ to be the associative algebra over
$\C((t))$ with generators $\beta_{t}(n), \gamma_{t}(n)\; (n\in \Z)$,
subject to relations
\begin{eqnarray*}
&&\beta_{t}(x_{1})\beta_{t}(x_{2})=\left(\frac{(q-1)t+qx_{2}-x_{1}}{(1-q)t+x_{2}-qx_{1}}\right)
\beta_{t}(x_{2})\beta_{t}(x_{1}),\\
&&\gamma_{t}(x_{1})\gamma_{t}(x_{2})=\left(\frac{(q-1)t+qx_{2}-x_{1}}{(1-q)t+x_{2}-qx_{1}}\right)
\gamma_{t}(x_{2})\gamma_{t}(v,x_{1}),\\
&&\beta_{t}(x_{1})\gamma_{t}(x_{2})-\left(\frac{(1-q)t+x_{2}-qx_{1}}{(q-1)t+qx_{2}-x_{1}}\right)
\gamma_{t}(x_{2})\beta_{t}(x_{1})
=x_{1}^{-1}\delta\left(\frac{x_{2}}{x_{1}}\right),
\end{eqnarray*}
where $\beta_{t}(x)=\sum_{n\in \Z}\beta_{t}(n)x^{-n-1}$,
$\gamma_{t}(x)=\sum_{n\in \Z}\gamma_{t}(n)x^{-n-1}$, and where when
$q\ne 1$, the rational-function coefficients are expanded in the
non-positive powers of $t$, e.g.,
$$\frac{1}{(1-q)t+x_{2}-qx_{1}}
=\sum_{i\ge 0}(1-q)^{-i-1}t^{-i-1}(x_{2}-qx_{1})^{i}\in
\C((t))[[x_{1},x_{2}]].$$ } \ed

By an {\em $A_{t,q}(\beta\gamma)$-module} we mean a $\C((t))$-module
$W$ on which $\beta_{t}(n),\gamma_{t}(n)$ for $n\in \Z$ act as
linear operators, satisfying the condition that for any $w\in W$,
$\beta_{t}(n)w=\gamma_{t}(n)w=0$ for $n$ sufficiently large and
those defining relations after applied to $w$ hold.
 A {\em vacuum
$A_{t,q}(\beta\gamma)$-module} is an $A_{t,q}(\beta\gamma)$-module
$W$ equipped with a vector $w_{0}$ that generates $W$ such that
$$\beta_{t}(n)w_{0}=\gamma_{t}(n)w_{0}=0\ \ \ \mbox{ for }n\ge 0.$$

Let $\tilde{A}$ be the free associative algebra over $\C((t))$ with
generators $\tilde{\beta}_{t}(n),\tilde{\gamma}_{t}(n)$ for $n\in
\Z$. Then an $A_{t,q}(\beta\gamma)$-module amounts to an
$\tilde{A}$-module $W$ such that for any $w\in W$,
$\tilde{\beta}_{t}(n)w=\tilde{\gamma}_{t}(n)w=0$ for $n$
sufficiently large and such that the three corresponding relations
after applied to $w$ hold. Define
$$\deg 1=0,\ \ \ \deg \tilde{\beta}_{t}(n)=
\deg \tilde{\gamma}_{t}(n)=-n-\frac{1}{2}\ \ \ \mbox{ for }n\in
\Z,$$ to make $\tilde{A}$ a $\frac{1}{2}\Z$-graded algebra, where
the degree-$k$ subspace is denoted by $\tilde{A}(k)$ for $k\in
\frac{1}{2}\Z$. We define an increasing filtration
${\mathcal{F}}=\{F_{k}\}_{k\in \frac{1}{2}\Z}$ of $\tilde{A}$ by
$F_{k}=\oplus_{p\le k} \tilde{A}(p)$ for $k\in \frac{1}{2}\Z$.
Clearly,
$$F_{p}\cdot F_{k}\subset F_{p+k}\ \ \ \mbox{ for }p,k\in \frac{1}{2}\Z.$$

\br{rvsalgebra} {\em Let $B$ be the associative algebra over $\C$
with generators $a_{n},b_{n}$ $(n\in \Z)$, subject to relations
$$a_{m}a_{n}=-a_{n}a_{m},\ \ \ \ b_{m}b_{n}=-b_{n}b_{m},\ \ \ \
a_{m}b_{n}+b_{n}a_{m}=\delta_{m+n+1,0}$$ for $m,n\in \Z$. In terms
of the generating functions
$$a(x)=\sum_{n\in \Z}a_{n}x^{-n-1},\ \ \ b(x)=\sum_{n\in
\Z}b_{n}x^{-n-1},$$ the above defining relations amount to
\begin{eqnarray*}
&&a(x_{1})a(x_{2})=-a(x_{2})a(x_{1}), \ \ \ \
b(x_{1})b(x_{2})=-b(x_{2})b(x_{1}),\\
&&\hspace{1cm}a(x_{1})b(x_{2})+b(x_{2})a(x_{1})
=x_{1}^{-1}\delta\left(\frac{x_{2}}{x_{1}}\right).
\end{eqnarray*}
Let $J$ be the left ideal of $B$, generated by $a_{n},b_{n}$ for
$n\ge 0$. Set
$$V_{B}=B/J,$$
a (left) $B$-module, set ${\bf 1}=1+J\in V_{B}$, and set
$$a=a_{-1}{\bf 1},\ \ \ \ b=b_{-1}{\bf 1}\in V_{B}.$$
It is well known (cf. \cite{ffr}) that $V_{B}$ is an irreducible
$B$-module. It follows that if $U$ is a nonzero $B$-module with a
vector $u_{0}$ satisfying the condition that $U=Bu_{0}$ and
$a_{n}u_{0}=b_{n}u_{0}=0$ for $n\ge 0$, then $U$ must be isomorphic
to $V_{B}$. It is also well known (see \cite{ffr}) that there exists
a vertex superalgebra structure on $V_{B}$, which is uniquely
determined by the conditions that ${\bf 1}$ is the vacuum vector and
that $Y(a,x)=a(x)$, $Y(b,x)=b(x)$.} \er

\bp{pirreducibility} Assume $q\ne 1$. Let $(W,w_{0})$ be a nonzero
vacuum $A_{t,q}(\beta\gamma)$-module. Then $F_{-1/2}w_{0}=0$ and $W$
is irreducible with $\End_{A_{t,q}(\beta\gamma)}(W)=\C((t))$. \ep

\begin{proof} From definition, $W$ is an $\tilde{A}$-module satisfying that
$\tilde{\beta}_{t}(n)w_{0}=\tilde{\gamma}_{t}(n)w_{0}=0$ for $n\ge
0$ and $W=\tilde{A}w_{0}$. We define an increasing sequence $W[k]$
with $k\in \frac{1}{2}\Z$ as follows: For $k<0$, set $W[k]=0$, and
for $k\in \frac{1}{2}\N$, let $W[k]$ be the span of the vectors
$$a^{(1)}(-m_{1})\cdots a^{(r)}(-m_{r})w_{0}$$
for $r\ge 0,\; a^{(1)},\dots,a^{(r)}\in \{
\tilde{\beta}_{t},\tilde{\gamma}_{t}\},\; m_{i}\ge 1$ with
$$\deg a^{(1)}(-m_{1})+\cdots +\deg a^{(r)}(-m_{r})=
(m_{1}-1/2)+\cdots +(m_{r}-1/2)\le k.$$ In the following we prove
that $W[k]=F_{k}w_{0}$ for $k\in \frac{1}{2}\Z$.

{}From definition, we have $W[0]=\C((t))w_{0}$ and
\begin{eqnarray}\label{e1234}
a(m)W[k]\subset W[k-m-1/2]\ \ \mbox{ for
 }a\in \{\tilde{\beta}_{t},\tilde{\gamma}_{t}\},\; m<0,\; k\in \frac{1}{2}\Z.
 \end{eqnarray}
Next, we show that this is also true for $m\ge 0$. Notice that for
$a,b\in \{ \tilde{\beta}_{t},\tilde{\gamma}_{t}\},\; m,n\in \Z,\;
w\in W$, from the defining relations in Definition \ref{dtbetagamma}
we have
$$a(m)b(n)w=-b(n)a(m)w+\sum_{i,j\ge 0,\; i+j\ge 1}f_{i,j}(t)
b(n+i)a(m+j)w+\lambda \delta_{m+n+1,0}w,$$ where $f_{i,j}(t)\in
\C((t))$, $\lambda=0$, or $1$. Then using induction on $k$, we can
show that (\ref{e1234}) also holds for $m\ge 0$, noting that
$a(m)w_{0}=0$ for $m\ge 0$. It follows that
$F_{k}w_{0}=F_{k}W[0]\subset W[k]$ for $k\in \frac{1}{2}\Z$. {}From
definition, we also have $W[k]\subset F_{k}w_{0}$. Therefore,
$F_{k}w_{0}=W[k]$. Consequently, $F_{k}w_{0}=W[k]=0$ for $k<0$. In
particular, we have $F_{-1/2}w_{0}=0$.

Now, we prove $\End_{A_{t,q}(\beta\gamma)}(W)=\C((t))$. We see that
the subspaces $W[k]$ $(k\in \frac{1}{2}\N)$ form an increasing
filtration of $W$, satisfying that $F_{p}W[k]\subset W[k+p]$ for
$k,p\in \frac{1}{2}\Z$. Form the associated graded space ${\rm
Gr}_{\mathcal{F}}(W)=\oplus_{k\in \frac{1}{2}\N}(W[k]/W[k-1/2])$,
which is naturally an $\tilde{A}$-module as ${\rm
Gr}_{\mathcal{F}}(\tilde{A})\simeq \tilde{A}$. Let $\rho:
\tilde{A}\rightarrow \End ({\rm Gr}_{\mathcal{F}}(W))$ be the
corresponding algebra homomorphism. On ${\rm Gr}_{\mathcal{F}}(W)$,
the following relations hold:
\begin{eqnarray*}
&&\rho(\tilde{\beta}_{t}(x))\rho(\tilde{\beta}_{t}(z))
=-\rho(\tilde{\beta}_{t}(z))\rho(\tilde{\beta}_{t}(x)),
\ \ \ \
\rho(\tilde{\gamma}_{t}(x))\rho(\tilde{\gamma}_{t}(z))
=-\rho(\tilde{\gamma}_{t}(z))\rho(\tilde{\gamma}_{t}(x)),\ \ \ \\
&&\hspace{2cm} \rho(\tilde{\beta}_{t}(x))\rho(\tilde{\gamma}_{t}(z))
+\rho(\tilde{\gamma}_{t}(z))\rho(\tilde{\beta}_{t}(x))
=x_{1}^{-1}\delta\left(\frac{x_{2}}{x_{1}}\right).
\end{eqnarray*}
We see that $\rho(\tilde{A})$ is a homomorphism image of
$\C((t))\otimes B$ (where $B$ is the algebra defined in Remark
\ref{rvsalgebra}), so that ${\rm Gr}_{\mathcal{F}}(W)$ is naturally
a $(\C((t))\otimes B)$-module. Since $W=\tilde{A}w_{0}$, we have
${\rm Gr}_{\mathcal{F}}(W)=\tilde{A}w_{0}$ with $w_{0}$ identified
with $w_{0}+W[-1/2]\in W[0]/W[-1/2]$. Then
$${\rm Gr}_{\mathcal{F}}(W)=(\C((t))\otimes B)w_{0}.$$
 {}From Remark \ref{rvsalgebra}, the $B$-submodule of
${\rm Gr}_{\mathcal{F}}(W)$, generated from $w_{0}$, is irreducible
and isomorphic to $V_{B}$. As $V_{B}$ is of countable dimension over
$\C$, we have $\End_{B}(V_{B})=\C$. Then one can show (cf.
\cite{li-simple}) that $\C((t))\otimes V_{B}$ is an irreducible
$\C((t))\otimes B$-module. It follows that ${\rm
Gr}_{\mathcal{F}}(W)\simeq \C((t))\otimes V_{B}$ as a
$\C((t))\otimes B$-module. Thus ${\rm Gr}_{\mathcal{F}}(W)$ is an
irreducible $\tilde{A}$-module.

Set
$$\Omega(W)=\{ w\in W\;|\; \beta_{t}(n)w=0=\gamma_{t}(n)w\ \ \ \mbox{
for }n\ge 0\}.$$  It is known that $\{v\in V_{B}\;|\;
a_{n}v=b_{n}v=0\ \ \mbox{ for }n\ge 0\}=\C{\bf 1}$. Then
$$\{w\in \C((t))\otimes V_{B}\;|\; a_{n}w=b_{n}w=0\ \ \mbox{ for }n\ge
0\}=\C((t)){\bf 1}.$$  Using this and the filtration $\mathcal{F}$
we obtain $\Omega(W)=\C((t))w_{0}$. Notice that for any endomorphism
$\psi$ of $W$, $\psi(w_{0})\in \Omega(W)$ and $\psi(w_{0})$
determines $\psi$ uniquely. Then it follows that
$\End_{A_{t,q}(\beta\gamma)}(W)=\C((t))$.

To prove that $W$ is irreducible, let $M$ be any submodule of $W$.
Then $M\cap W[k]$ with $k\in \frac{1}{2}\Z$ form an increasing
filtration of $M$ and the associated graded space ${\rm
Gr}_{\mathcal{F}}(M)$ can be considered canonically as a subspace of
${\rm Gr}_{\mathcal{F}}(W)$. It is clear that ${\rm
Gr}_{\mathcal{F}}(M)$ is an $A$-submodule. As ${\rm
Gr}_{\mathcal{F}}(W)$ is an irreducible $A$-module, we must have
either ${\rm Gr}_{\mathcal{F}}(M)=0$ or ${\rm
Gr}_{\mathcal{F}}(M)={\rm Gr}_{\mathcal{F}}(W)$. If ${\rm
Gr}_{\mathcal{F}}(M)=0$, we have $M\cap W[k]=M\cap K[k-1/2]$ for all
$k\in \frac{1}{2}\Z$. Since $W[k]=0$ for $k$ sufficiently negative,
we have $M\cap W[k]=0$ for all $k$. Thus $M=0$. On the other hand,
if ${\rm Gr}_{\mathcal{F}}(M)={\rm Gr}_{\mathcal{F}}(W)$, we have
$M\cap W[k]+W[k-1/2]=W[k]$ for all $k$. Using induction we get
$W[k]\subset M$ for all $k$. Thus $M=W$. This proves that $W$ is
irreducible, concluding the proof.
\end{proof}

The following gives the existence of a nonzero vacuum
$V_{t,q}(\beta\gamma)$-module:

\bp{pnonzero} Let $V_{B}$ be the vertex superalgebra as in Remark
\ref{rvsalgebra}. There exists linear maps
$$\Phi^{\pm}(t): V_{B}\rightarrow V_{B}\otimes \C((t))$$ satisfying the condition
that
\begin{eqnarray*}
&&\Phi^{\pm}(t){\bf 1}={\bf 1}, \ \ \ \Phi^{\pm}(t)(a)=a\otimes
t^{\pm 1},\ \ \
\Phi^{\pm}(t)(b)=b\otimes t^{\mp 1},\\
&&\Phi^{\pm}(x_{1})Y(v,x_{2})=Y(\Phi^{\pm}(x_{1}-x_{2})v,x_{2})\Phi^{\pm}(x_{1})
\ \ \ \mbox{ for }v\in V_{B},\\
&&\Phi^{\pm}(x_{1})\Phi^{\pm}(x_{2})=\Phi^{\pm}(x_{2})\Phi^{\pm}(x_{1}),\\
&&\Phi^{+}(x)\Phi^{-}(x)=\Phi^{-}(x)\Phi^{+}(x)=1.
\end{eqnarray*}
Furthermore, if $q\ne 1$, the assignment
\begin{eqnarray*}
&&\beta_{t}(x)=(1-q)(t+x)Y(a,qx)\Phi((1-q)t+x),\\
&&\gamma_{t}(x)=Y(b,qx)\Phi^{-1}((1-q)t+x)
\end{eqnarray*}
defines a vacuum $A_{t,q}(\beta\gamma)$-module structure on
$V_{B}\otimes \C((t))$. \ep

\begin{proof} It is similar to the proof of a similar result in Section 4
of \cite{li-qva2}. Equip $\C((t))$ with the vertex algebra structure with $1$ as the
vacuum vector and with
$$Y(f(t),x)g(t)=(e^{-x(d/dt)}f(t))g(t)=f(t-x)g(t)
\ \ \ \mbox{ for }f(t),g(t)\in \C((t)).$$ Furthermore, equip
$V_{B}\otimes \C((t))$ with the vertex superalgebra structure by
tensor product over $\C$. We have
$$Y(a\otimes t^{\pm 1},x)=Y(a,x)\otimes (t+x)^{\pm 1},\ \
Y(b\otimes t^{\mp 1},x)=Y(b,x)\otimes (t+x)^{\mp 1}.$$ It is
straightforward to show that the assignments
$$a(x)\mapsto Y(a\otimes t^{\pm 1},x), \ \ b(x)\mapsto Y(b\otimes t^{\mp
1},x)$$ give two $B$-module structures on $V_{B}\otimes \C((t))$. It
follows from the universal property of $V_{B}$ that there exist
$B$-module homomorphisms $\Phi^{\pm}: V_{B}\rightarrow V_{B}\otimes
\C((t))$ such that $\Phi^{\pm}({\bf 1})={\bf 1}\otimes 1$. Since
$a,b$ generate $V_{B}$ as a vertex algebra, it follows that
$\Phi^{\pm}$ are vertex algebra homomorphisms. We have
$$ \Phi^{\pm}(a)=\Phi^{\pm}(a_{-1}{\bf 1})=\Res_{x}x^{-1}(Y(a,x)\otimes
(t+x)^{\pm 1})({\bf 1}\otimes 1) =a\otimes t^{\pm 1},\ \ \
\Phi^{\pm} (b)=b\otimes t^{\mp 1}.$$ Write $\Phi^{\pm}$ as $\Phi(t)$
and $\Phi^{-}(t)$, indicating the dependence on $t$. Then
$\Phi^{\pm}(t)$ meet all the requirements.

As for the last assertion, note that $\Phi((1-q)t+x)$ makes sense as
$\Phi(x)(v)\in V\otimes \C((x))$. We have
\begin{eqnarray*}
& &Y(a,qx_{1})\Phi((1-q)t+x_{1})Y(a,qx_{2})\Phi((1-q)t+x_{2})\\
&=&((1-q)t+x_{1}-qx_{2})Y(a,qx_{1})Y(a,qx_{2})\Phi((1-q)t+x_{1})\Phi((1-q)t+x_{2})\\
&=&((q-1)t+qx_{2}-x_{1})Y(a,qx_{2})Y(a,qx_{1})\Phi((1-q)t+x_{1})\Phi((1-q)t+x_{2})\\
&=&\left(\frac{(q-1)t+qx_{2}-x_{1}}{(1-q)t+x_{2}-qx_{1}}\right)
Y(a,qx_{2})\Phi((1-q)t+x_{2})Y(a,qx_{1})\Phi((1-q)t+x_{1}),
\end{eqnarray*}
\begin{eqnarray*}
& &Y(b,qx_{1})\Phi^{-}((1-q)t+x_{1})Y(b,qx_{2})\Phi^{-}((1-q)t+x_{2})\\
&=&((1-q)t+x_{1}-qx_{2})Y(b,qx_{1})Y(b,qx_{2})\Phi^{-}((1-q)t+x_{1})\Phi^{-}((1-q)t+x_{2})\\
&=&((q-1)t+qx_{2}-x_{1})Y(b,qx_{2})Y(b,qx_{1})\Phi^{-}((1-q)t+x_{1})\Phi^{-}((1-q)t+x_{2})\\
&=&\left(\frac{(q-1)t+qx_{2}-x_{1}}{(1-q)t+x_{2}-qx_{1}}\right)
Y(b,qx_{2})\Phi^{-}((1-q)t+x_{2})Y(b,qx_{1})\Phi^{-}((1-q)t+x_{1}),
\end{eqnarray*}
\begin{eqnarray*}
& &Y(a,qx_{1})\Phi((1-q)t+x_{1})Y(b,qx_{2})\Phi^{-}((1-q)t+x_{2})\\
&&\ \ -\frac{(1-q)t+x_{2}-qx_{1}}{(q-1)t+qx_{2}-x_{1}}Y(b,x_{2})
\Phi^{-}((1-q)t+x_{2})Y(a,qx_{1})\Phi((1-q)t+x_{1})\\
&=&((1-q)t+x_{1}-qx_{2})^{-1}Y(a,qx_{1})Y(b,qx_{2})\Phi((1-q)t+x_{1})\Phi^{-}((1-q)t+x_{2})\\
& &-((q-1)t+qx_{2}-x_{1})^{-1}Y(b,qx_{2})Y(a,qx_{1})
\Phi^{-}((1-q)t+x_{2})\Phi((1-q)t+x_{1})\\
&=&
((1-q)t+x_{1}-qx_{2})^{-1}x_{2}^{-1}\delta\left(\frac{x_{1}}{x_{2}}\right)
\Phi^{-}((1-q)t+x_{2})\Phi((1-q)t+x_{1})\\
&=&(1-q)^{-1}(t+x_{1})^{-1}x_{2}^{-1}\delta\left(\frac{x_{1}}{x_{2}}\right).
\end{eqnarray*}
This proves that $V_{B}\otimes \C((t))$ is an
$A_{t,q}(\beta\gamma)$-module with the given action. Let $M$ be the
$A_{t,q}(\beta\gamma)$-submodule of $V_{B}\otimes \C((t))$,
generated from ${\bf 1}\otimes 1$. It is clear that $M$ is a vacuum
$A_{t,q}(\beta\gamma)$-module. Now, it suffices to prove that
$M=V_{B}\otimes \C((t))$. As $V_{B}$ is an irreducible $B$-module,
we have $V_{B}=B\cdot{\bf 1}$, so that $$V_{B}\otimes
\C((t))=(B\otimes \C((t)))({\bf 1}\otimes 1).$$ Then it suffices to
prove that $M$ is stable under the action of $B$. With
$\Phi^{\pm}(x){\bf 1}={\bf 1}$, using the commutation relations (and
induction), we see that $M$ is stable under the actions of
$\Phi^{\pm}((1-q)t+x)$. Note that by definition $M$ is stable under
the actions of $\beta_{t}(x)$ and $\gamma_{t}(x)$. Consequently, $M$
is stable under the actions of $Y(a,x)$ and $Y(b,x)$. Thus $M$ is
stable under the action of $B$. Therefore, we have $M=V_{B}\otimes
\C((t))$, proving that $V_{B}\otimes \C((t))$ is a vacuum
$A_{t,q}(\beta\gamma)$-module.
\end{proof}

We now construct a universal vacuum $A_{t,q}(\beta\gamma)$-module.
First, set $\tilde{J}=\tilde{A}F_{-1/2}$ (recall Proposition
\ref{pirreducibility}), a left ideal of $\tilde{A}$.
 Then consider the
quotient $\tilde{A}/\tilde{J}$, a (left) $\tilde{A}$-module. One
sees that for any $w\in \tilde{A}/\tilde{J}$,
$\tilde{\beta}_{t}(n)w=\tilde{\gamma}_{t}(n)w=0$ for $n$
sufficiently large, as for any $a\in \tilde{A}$,
$\tilde{\beta}_{t}(n)a, \; \tilde{\gamma}_{t}(n)a\in F_{-1/2}$ for
$n$ sufficiently large.

\bd{duniversal} {\em Let $V_{t,q}(\beta\gamma)$ be the quotient of
$\tilde{A}/\tilde{J}$ modulo the relations corresponding to the
defining relations of $A_{t,q}(\beta\gamma)$. We set ${\bf
1}=1+\tilde{J}\in V_{t,q}(\beta\gamma)$.} \ed

{}From the construction, $(V_{t,q}(\beta\gamma),{\bf 1})$ is a
vacuum $A_{t,q}(\beta\gamma)$-module and it is universal in the
obvious sense.  It then follows from Propositions
\ref{pirreducibility} and \ref{pnonzero} that $V_{t,q}(\beta\gamma)$
is irreducible (nonzero) and that every nonzero vacuum module is
isomorphic to $V_{t,q}(\beta\gamma)$.

Now we are ready to present the main result of this section:

\bt{tqvabetagamma} Assume $q\ne 1$. There exists a weak quantum
vertex $\C((t))$-algebra structure on $V_{t,q}(\beta\gamma)$ with
${\bf 1}$ as the vacuum vector and with
$$Y(\beta_{t}(-1){\bf 1},x)=\beta_{t}(x),\ \ \ \
Y(\gamma_{t}(-1){\bf 1},x)=\gamma_{t}(x).$$ Furthermore, such a weak
quantum vertex $\C((t))$-algebra structure is unique and
non-degenerate. \et

\begin{proof} We shall apply Theorem \ref{tconstruction-refine}. Set
$$U=\C((t))\beta_{t}+\C((t))\gamma_{t}\subset V_{t,q}(\beta\gamma)$$
and define
$$Y_{0}(f(t)\beta_{t},x)=f(t+x)\beta_{t}(x),\ \
Y_{0}(f(t)\gamma_{t},x)=f(t+x)\gamma_{t}(x) \ \ \mbox{ for }f(t)\in
\C((t)).$$ Set $U(x)=\{ Y_{0}(u,x)\;|\; u\in U\}$. It is clear that
$U(x)$ is $\S_{t}(x_{1},x_{2})$-local, so $U(x)$ generates a
nonlocal vertex algebra $\<U(x)\>$ over $\C((t))$ inside
$\E(V_{t,q}(\beta\gamma))$. Furthermore, $\C((t))((x))\<U(x)\>$ is a
nonlocal vertex algebra over $\C((t))$. By Proposition
\ref{prelation-am} we have
\begin{eqnarray*}
&&Y_{\E}(\beta_{t}(x),x_{1})Y_{\E}(\beta_{t}(x),x_{2})
=\left(\frac{(q-1)(t+x)+qx_{2}-x_{1}}{(1-q)(t+x)+x_{2}-qx_{1}}\right)
Y_{\E}(\beta_{t}(x),x_{2})Y_{\E}(\beta_{t}(x),x_{1}),\\
&&Y_{\E}(\gamma_{t}(x),x_{1})Y_{\E}(\gamma_{t}(x),x_{2})
=\left(\frac{(q-1)(t+x)+qx_{2}-x_{1}}{(1-q)(t+x)+x_{2}-qx_{1}}\right)
Y_{\E}(\gamma_{t}(x),x_{2})Y_{\E}(\gamma_{t}(x),x_{1}),\\
&&Y_{\E}(\beta_{t}(x),x_{1})Y_{\E}(\gamma_{t}(x),x_{2})
-\left(\frac{(1-q)(t+x)+x_{2}-qx_{1}}{(q-1)(t+x)+qx_{2}-x_{1}}\right)
Y_{\E}(\gamma_{t}(x),x_{2})Y_{\E}(\beta_{t}(x),x_{1})\\
&&\hspace{1cm}=x_{1}^{-1}\delta\left(\frac{x_{2}}{x_{1}}\right).
\end{eqnarray*}
Define a $\C((t_{1}))$-module structure on $\C((t))((x))\<U(x)\>$
with $f(t_{1})\in \C((t_{1}))$ acting as $f(t+x)$. Then
$\C((t))((x))\<U(x)\>$ is an $A_{t_{1},q}(\beta\gamma)$-module with
$\beta_{t}(z)$ and $\gamma_{t}(z)$ acting as
$Y_{\E}(\beta_{t}(x),z)$ and $Y_{\E}(\gamma_{t}(x),z)$,
respectively. Furthermore, the submodule generated from
$1_{V_{t,q}(\beta\gamma)}$ is a vacuum
$A_{t_{1},q}(\beta\gamma)$-module. It follows that there is a
$\C$-linear map $\psi$ from $V_{t,q}(\beta\gamma)$ to
$\C((t))((x))\<U(x)\>$ such that
\begin{eqnarray*}
&&\ \ \ \ \psi({\bf 1})=1_{V},\ \ \ \ \psi(f(t)v)=f(t_{1})\psi(v), \\
&&\psi( \beta_{t}(z)v)=Y_{\E}(\beta_{t}(x),z)\psi(v),\ \ \psi(
\gamma_{t}(z)v)=Y_{\E}(\gamma_{t}(x),z)\psi(v)
\end{eqnarray*}
for $f(t)\in \C((t)),\;v\in V_{t,q}(\beta\gamma)$.  Now the first
assertion follows from Theorem \ref{tconstruction-refine}.

Next, we show that $V_{t,q}(\beta\gamma)$ is non-degenerate by using
Proposition \ref{pnon-degenerate}. Recall the $\frac{1}{2}\Z$-graded
free algebra $\tilde{A}$ over $\C((t))$. By Proposition
\ref{pirreducibility}, $V_{t,q}(\beta\gamma)$ is an irreducible
$\tilde{A}$-module with $\End_{\tilde{A}}
(V_{t,q}(\beta\gamma))=\C((t))$. As $\beta_{t},\gamma_{t}$ generate
$V_{t,q}(\beta\gamma)$, it follows that $V_{t,q}(\beta\gamma)$ as a
$V_{t,q}(\beta\gamma)$-module is irreducible with
$\End_{V_{t,q}(\beta\gamma)} (V_{t,q}(\beta\gamma))=\C((t))$. Now,
by Proposition \ref{pnon-degenerate}, $V_{t,q}(\beta\gamma)$ is
non-degenerate.
\end{proof}

Regarding the relationship between $A_{q}(\beta\gamma)$-modules and
the quantum vertex $\C((t))$-algebra $V_{t,q}(\beta\gamma)$ we have:

\bp{pbetagamma-module} Assume $q\ne 1$. Let $W$ be an
$A_{q}(\beta\gamma)$-module. There exists a type zero module
structure for the quantum vertex $\C((t))$-algebra
$V_{t,q}(\beta\gamma)$ with
$$Y_{W}(\beta_{t}(-1){\bf 1},x)=\beta(x),\ \ \ \
Y_{W}(\gamma_{t}(-1){\bf 1},x)=\gamma(x).$$
 \ep

\begin{proof} From the defining relations of $A_{q}(\beta\gamma)$,
one sees that the generating functions $\beta(x)$ and $\gamma(x)$
form an $\S(x_{1},x_{2})$-local subset of $\E(W)$. Thus by Theorem
\ref{tquasi-main}, $\{\beta(x),\gamma(x)\}$ generates a nonlocal
vertex algebra $K$ over $\C$ with $W$ as a module. Furthermore, by
Theorem \ref{tconcrete-general}, $\C((x))K$ is a weak quantum vertex
$\C((t))$-algebra with $W$ as a type zero module, where $f(t)\in
\C((t))$ acts on $\C((x))K$ as $f(x)$. In view of Proposition
\ref{prelation-am} we have
\begin{eqnarray*}
&&Y_{\E}(\beta(z),x_{1})Y_{\E}(\beta(z),x_{2})
=\left(\frac{(1-q)t+x_{1}-qx_{2}}{(q-1)t+qx_{1}-x_{2}}\right)
Y_{\E}(\beta(z),x_{2})Y_{\E}(\beta(z),x_{1}),\\
&&Y_{\E}(\gamma(z),x_{1})Y_{\E}(\gamma(z),x_{2})
=\left(\frac{(1-q)t+x_{1}-qx_{2}}{(q-1)t+qx_{1}-x_{2}}\right)
Y_{\E}(\gamma(z),x_{2})Y_{\E}(\gamma(z),x_{1}),\\
&&Y_{\E}(\beta(z),x_{1})Y_{\E}(\gamma(z),x_{2})
-\left(\frac{(q-1)t+qx_{1}-x_{2}}{(1-q)t+x_{1}-qx_{2}}\right)
Y_{\E}(\gamma(z),x_{2})Y_{\E}(\beta(z),x_{1})\\
&&\hspace{2cm}=x_{1}^{-1}\delta\left(\frac{x_{2}}{x_{1}}\right).
\end{eqnarray*}
Thus, $\C((x))K$ is an $A_{t,q}(\beta\gamma)$-module with $f(t)\in
\C((t))$ acting as $f(x)$ and with $\beta_{t}(n)$, $\gamma_{t}(n)$
acting as $\beta(z)_{n},\gamma(z)_{n}$ for $n\in \Z$, respectively.
We also have $\beta(z)_{n}1_{W}=0=\gamma(z)_{n}1_{W}$ for $n\ge 0$.
It follows that there exists an $A_{t,q}(\beta\gamma)$-module
homomorphism $\pi$ {}from $V_{t,q}(\beta\gamma)$ to $\C((x))K$,
sending ${\bf 1}$ to $1_{W}$. That is, $\pi$ is a $\C((t))$-module
homomorphism satisfying the condition that $\pi({\bf 1})=1_{W}$,
$$\pi(Y(\beta_{t},z)v)=Y_{\E}(\beta(x),z)\pi(v), \ \ \
\pi(Y(\gamma_{t},z)v)=Y_{\E}(\gamma(x),z)\pi(v)  $$ for $v\in
V_{t,q}(\beta\gamma)$. It follows that $\pi$ is a homomorphism of
nonlocal vertex $\C((t))$-algebras. Consequently, $W$ is a module of
type zero for $V_{t,q}(\beta\gamma)$.
\end{proof}

\section*{Appendix}
In this Appendix we present some technical results on iota-maps,
which we use in the main body of this paper.

\bl{lfirst-form} For any $f(x_{1},x_{2})\in \F_{*}(x_{1},x_{2})$, we
have
\begin{eqnarray*}
&&\iota_{t_{1},x_{2}}\left(\iota_{t,x_{0}}f(t+x_{0},t)\right)|_{t=t_{1}+x_{2}}
=\iota_{t_{1},x_{2},x_{0}}f(t_{1}+x_{2}+x_{0},t_{1}+x_{2}),\\
&&\iota_{x_{1},x_{0}}\left(\iota_{t,x_{2},x_{0}}f(t+x_{2}+x_{0},t+x_{2})\right)|_{x_{2}=x_{1}-x_{0}}
=\iota_{t, x_{1},x_{0}}f(t+x_{1},t+x_{1}-x_{0}).\ \ \ \ \
\end{eqnarray*}
\el

\begin{proof}
If $f(x_{1},x_{2})\in \F[[x_{1},x_{2}]]$, it is clear as iota-maps
leave nonnegative power series unchanged. Now, we consider the case
with $f=1/p$ for $p(x_{1},x_{2})\in \F[x_{1},x_{2}]$. We have
$\iota_{t,x_{0}}(1/p(t+x_{0},t))=p(t+x_{0},t)^{-1}$, the inverse of
$p(t+x_{0},t)$ in $\F((t))((x_{0}))$. The substitution
$t=t_{1}+x_{2}$ is an algebra homomorphism from $\F((t))((x_{0}))$
into $\F((t_{1}))((x_{2}))((x_{0}))$. Thus
$\left(\iota_{t,x_{0}}(1/p(t+x_{0},t))\right)|_{t=t_{1}+x_{2}}$ is
the inverse of polynomial $p(t_{1}+x_{2}+x_{0},t_{1}+x_{2})$ in
 $\F((t_{1}))((x_{2}))((x_{0}))$.
On the other hand, we know that
$\iota_{t_{1},x_{2},x_{0}}(1/p(t_{1}+x_{2}+x_{0},t_{1}+x_{2}))$ is
also the inverse of $p(t_{1}+x_{2}+x_{0},t_{1}+x_{2})$  in
 $\F((t_{1}))((x_{2}))((x_{0}))$. This proves the first assertion.
The second assertion can be proved similarly.
\end{proof}

\bl{lsubcomm} a) Let $q(x_{1},x_{2})\in \F[x_{1},x_{2}]$ be such
that $q(x_{1},x_{1})\ne 0$. Then
\begin{eqnarray*}
\iota_{t,x_{1},x_{2}}(1/q(t+x_{1},t+x_{2}))
=\iota_{t,x_{2},x_{1}}(1/q(t+x_{1},t+x_{2}))\in
\F((t))[[x_{1},x_{2}]].
\end{eqnarray*}
b) For any $f(x_{1},x_{2})\in \F_{*}(x_{1},x_{2})$, we have
\begin{eqnarray*}
\iota_{t,x_{2},x_{1}}f(t+x_{1},t+x_{2})\in
\F((t))((x_{2}))[[x_{1}]],
\end{eqnarray*}
and there exists $k\in \N$ such that
$$\iota_{t,x_{1},x_{2}}(x_{1}-x_{2})^{k}f(t+x_{1},t+x_{2})
=\iota_{t,x_{2},x_{1}}(x_{1}-x_{2})^{k}f(t+x_{1},t+x_{2})$$ with
both sides lying in $\F((t))[[x_{1},x_{2}]]$, and such that
\begin{eqnarray*}
\iota_{t,x_{1},x_{0}}(x_{0}^{k}f(t+x_{1},t+x_{1}-x_{0}))
=\left(\iota_{t,x_{2},x_{1}}(x_{1}-x_{2})^{k}f(t+x_{1},t+x_{2})\right)|_{x_{2}=x_{1}-x_{0}}.
\end{eqnarray*}
 \el

\begin{proof} We have
$$q(t+x_{1},t+x_{2})=q(t,t)-x_{1}g(t,x_{1},x_{2})-x_{2}h(t,x_{1},x_{2})$$
for some $g,h\in \F[t,x_{1},x_{2}]$ where $q(t,t)\ne 0$. Then
\begin{eqnarray*}
&&\iota_{t,x_{1},x_{2}}(1/q(t+x_{1},t+x_{2}))\\
&=&\sum_{j\ge
0}\iota_{t,0}(1/q(t,t))^{j+1}(x_{1}g(t,x_{1},x_{2})+x_{2}h(t,x_{1},x_{2}))^{j}\\
&=&\iota_{t,x_{2},x_{1}}(1/q(t+x_{1},t+x_{2})),
\end{eqnarray*}
proving the first assertion. Let $f=g/p$ with $g\in
\F[[x_{1},x_{2}]]$, $p(x_{1},x_{2})\in \F[x_{1},x_{2}]$ (nonzero).
We have $p(x_{1},x_{2})=(x_{1}-x_{2})^{k}q(x_{1},x_{2})$ for some
$k\in \N,\; q(x_{1},x_{2})\in \F[x_{1},x_{2}]$ with
$q(x_{2},x_{2})\ne 0$. Then the second and the third assertions
follow immediately. As for the last assertion, we have
\begin{eqnarray*}
\iota_{t,x_{1},x_{0}}(1/q(t+x_{1},t+x_{1}-x_{0}))
=\left(\iota_{t,x_{2},x_{1}}(1/q(t+x_{1},t+x_{2})\right)|_{x_{2}=x_{1}-x_{0}},
\end{eqnarray*}
because both sides are the inverse of $q(t+x_{1},t+x_{1}-x_{0})$ in
$\F((t))((x_{1}))((x_{0}))$. Then the last assertion follows.
\end{proof}

\end{document}